%% file: homoclopoline_eng.tex
\documentclass[a4paper, 12pt]{amsart}
\usepackage[english]{babel}
\usepackage{amsfonts,amsmath,amssymb,amsthm,graphicx,afterpage,url}
\input{xy}
\xyoption{all}
\usepackage{epsfig}
\usepackage{overpic}
\usepackage{pgf,tikz}
\usetikzlibrary{arrows}
\usepackage{float}

\newcommand{\cycles}[1]{}
\newcommand{\algor}[1]{}
\newcommand{\invadraw}[1]{}

\newcommand{\compile}[1]{#1}

\graphicspath{{figs/}}

\def\R{{\mathbb R}} \def\Z{{\mathbb Z}}

\long\def\comment#1\endcomment{}

\theoremstyle{theorem}
\newtheorem{theorem}{Theorem}[section]

    \newtheorem{proposition}[theorem]{Proposition}

\theoremstyle{definition}

\newtheoremstyle{mydefinition}
  {3pt}
  {3pt}
  {\normalfont}
  {\parindent}
  {\bfseries}
  {.}
  { }
  {}

\theoremstyle{mydefinition}
\newtheorem{pr}[theorem]{Problem}

\makeatletter

\@addtoreset{figure}{section}
\makeatother

\usepackage{hyperref}
\hypersetup{
   colorlinks   = true, 
   urlcolor     = blue, 
   linkcolor    = blue, 
   citecolor   = red 
}

\begin{document}

\newpage
\title{Homotopy classification of closed polygonal lines: \\ results and problems}

\author{E. Alkin, O. Nikitenko and A. Skopenkov}

\thanks{\emph{E. Alkin, A. Skopenkov:} Moscow Institute of Physics and Technology.
\newline
\emph{O. Nikitenko:} Altay Technical University (Barnaul).
\newline
\emph{A. Skopenkov:} Independent University of Moscow, \url{https://users.mccme.ru/skopenko}.
\newline
We are grateful to A. Miroshnikov, Yu. Khromin for useful discussions, and to A. Miroshnikov for the preparation of some figures.}

\date{}



\begin{abstract}
In this text we expose (as a sequence of problems) basic cases of some fundamental ideas and methods of mathematics. 
Namely, of homotopy, degree, fundamental group, covering, Whitehead invariant, etc.
This is done by considering the elementary example: closed polygonal lines in a subset of the plane. 
Although these ideas and methods are parts of topology, they are used in many other areas including computer science.
 
This text is expository and is accessible to mathematicians not specialized in the area (and to students).
\end{abstract}

\maketitle
\tableofcontents

\input{obstr-hoclpoli_eng}

{\it Books, surveys, and expository papers in this list are marked by the stars.}

\end{document}


\title{Гомотопическая классификация замкнутых ломаных} 
 
\author{Э. Алкин и А. Скопенков}

\address{\emph{Э. Алкин, А. Скопенков:} Московский Физико-Технический Институт.
\newline
\emph{А. Скопенков:} Независимый Московский Университет, \url{https://users.mccme.ru/skopenko}.}

\thanks{Благодарим А. Мирошникова, О. Никитенко и Ю. Хромина за полезные обсуждения, а также А. Мирошникова за подготовку некоторых рисунков.} 
   
\date{}

\begin{abstract}
Мы излагаем 
базовые случаи некоторых фундаментальных идей и методов математики (гомотопии, степени, фундаментальной группы, накрытия, инварианта Уайтхеда и т.д.). 
Это делается на основе элементарного примера: замкнутые ломаные в подмножестве плоскости.
Хотя эти идеи и методы являются частью топологии, они используются во многих других областях, включая информатику. 
\end{abstract}

\maketitle
\tableofcontents


\section{Введение}\label{s:intr}

Неформально говоря, две замкнутые ломаные в~подмножестве плоскости называются \emph{гомотопными}, если одна может быть преобразована в~другую <<непрерывной деформацией>> в этом подмножестве.
Это понятие возникло, в частности, при построении непрерывных ветвей функции $\sqrt z$ комплексного аргумента на плоскости без точки (и других многозначных аналитических функций на других подмножествах плоскости): если две ломаные гомотопны на плоскости без точки, то построенные по ним две ветви дают одинаковые разности значений в начальной и конечной точке. 
(Дальнейший текст не зависит от этого предложения, поэтому, если оно Вам непонятно, читайте дальше.)

Приведем строгие определения.

\textbf{Ориентированной циклической последовательностью} называется упорядоченный набор с точностью до циклического сдвига.
\textbf{Замкнутой} ориентированной \textbf{ломаной} на плоскости называется\footnote{Это определение не согласовано с определением замкнутой ломаной в статье \cite{AMS'}, но мы считаем его более удобным в рамках данного материала.} ориентированная циклическая последовательность точек плоскости (не обязательно различных).
Далее слова <<ориентированная>> и <<на плоскости>> пропускаются.\footnote{Таким образом, замкнутая ломаная (определенная здесь) не является подмножеством плоскости.
Тем не менее, иногда мы работаем с замкнутой ломаной $A_1\ldots A_m$ как с объединением отрезков $A_iA_{i+1}$, например, мы пишем <<ломаная, не проходящая через точку>>.}

\begin{figure}[h]\centering
\compile{
\includegraphics[scale=0.5]{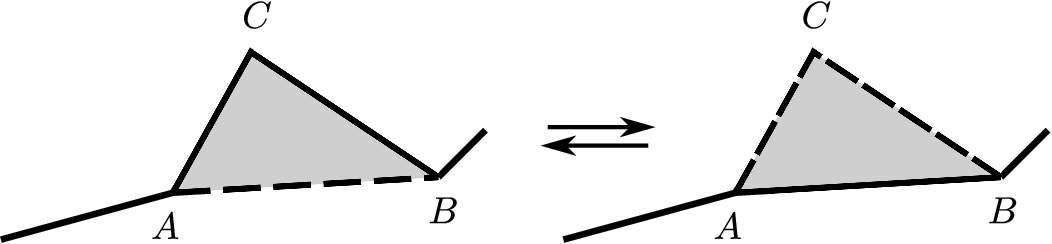}
}
\caption{Элементарное сокращение}
\label{f:elem}
\end{figure}

Пусть $N$ "--- подмножество плоскости.
\emph{Элементарным сокращением} в~$N$ замкнутой ломаной в~$N$ называется удаление такой ее вершины~$B$, что для соседних с~$B$ вершин $A,C$ (которые не обязательно различны) выпуклая оболочка вершин $A,B,C$ \emph{содержится в~$N$} (см. рисунок~\ref{f:elem}).
\emph{Гомотопией в $N$} (кусочно-линейной) называется конечная последовательность замкнутых ломаных, для любых двух последовательных ломаных которой одна получена из другой элементарным сокращением в~$N$.
Две замкнутые ломаные называются (кусочно-линейно) \textbf{гомотопными в~$N$}, если существует гомотопия в $N$, первая и последняя ломаные в которой совпадают с заданными.

\emph{Далее, когда мы пишем о ломаной в некотором подмножестве плоскости, мы рассматриваем гомотопность (и элементарные сокращения) в том же подмножестве.}

\begin{example}\label{p:veplho-poly0b} Любые две одноточечные замкнутые ломаные на плоскости без точки гомотопны.
Действительно, пусть даны точки~$A$ и~$B$.
Возьмем точку~$C$, не лежащую ни на одной из прямых, соединяющих $A$ и $B$ с выколотой точкой.
Тогда нужная гомотопия "--- $A,AC,C,CB,B$.
\end{example}

Гомотопическая классификация замкнутых ломаных \emph{на плоскости} тривиальна.

\begin{proposition}\label{p:veplho-poly0a} Любые две замкнутые ломаные на плоскости го\-мо\-топны.
\end{proposition}

\begin{figure}[h]\centering
\compile{
\includegraphics[scale=0.6]{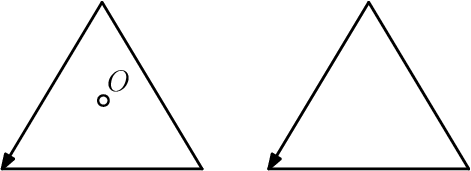}}
\caption{Две замкнутые ломаные на плоскости без точки $O$}
\label{f:guess}
\end{figure}

\begin{example}[см. рисунок~{\ref{f:guess}}]\label{p:veplho-poly1} 
Существуют негомотопные замкнутые ломаные на плоскости без точки.
\end{example}

Классификация ломаных с точностью до гомотопности на плоскости без точки приведена в \S\ref{s:minus}.

\begin{figure}[h]\centering
\compile{
\includegraphics[scale=0.8]{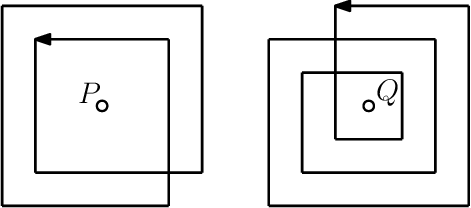} \qquad
\includegraphics[scale=0.8]{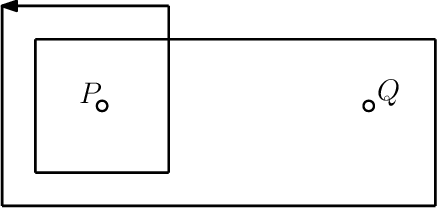}\\
\vspace{0.5cm}
\includegraphics[scale=0.9]{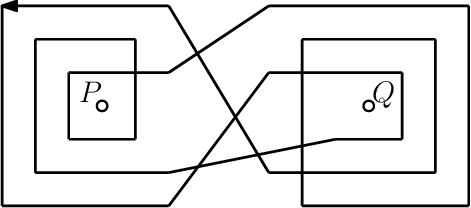}
}
\caption{Четыре замкнутые ломаные на плоскости без двух точек $P, Q$}
\label{f:2points-a2-b3}
\end{figure}

\begin{figure}[h]\centering
\compile{
\includegraphics[scale=0.25]{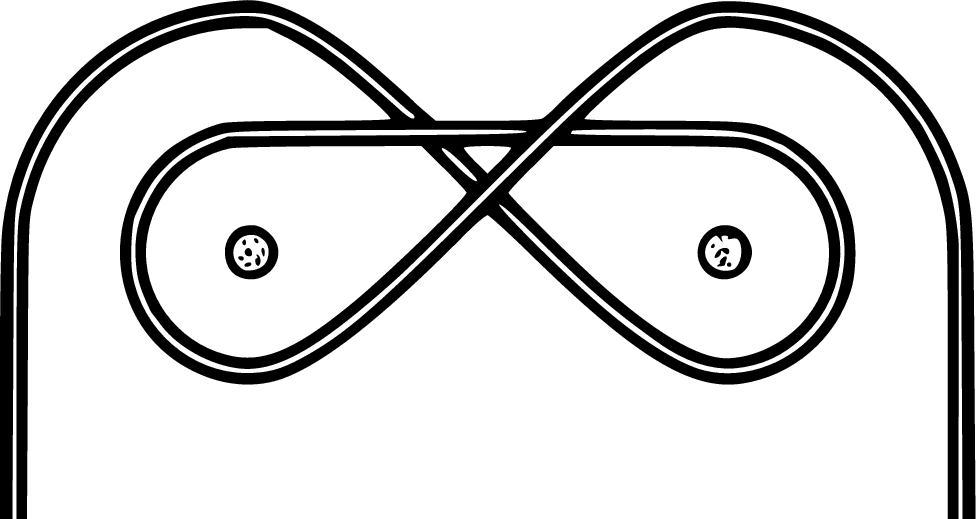}
\qquad
\qquad
\includegraphics[scale=0.6]{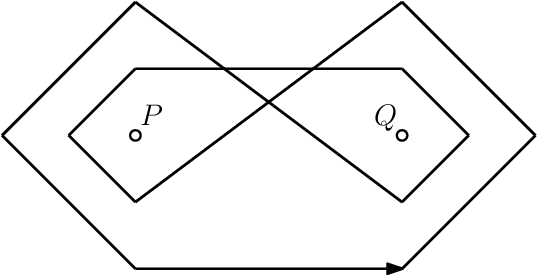}
}
\caption{Пример~\ref{veplho-borr}.a: веревка на двух гвоздях (слева); замкнутая ломаная на плоскости без двух точек $P,Q$ (справа)}
\label{f:poincare}
\end{figure}

\begin{figure}
\centering
\compile{
\includegraphics[scale=0.9]{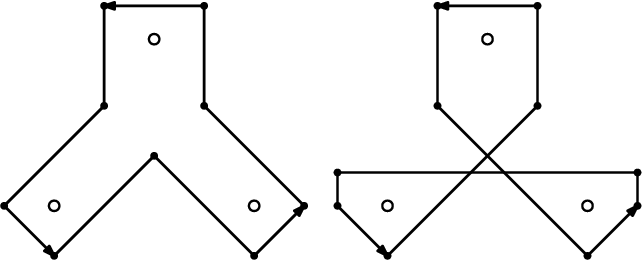}
}
\caption{Две замкнутые ломаные на плоскости без трех точек}
\label{f:3poin}
\end{figure}

Попробуйте придумать инварианты, которые различают замкнутые ломаные в 
примерах~\ref{p:veplho-poly1}--\ref{p:3points} и в задаче~\ref{3space} (в частности, на рисунках~\ref{f:guess}--\ref{f:poincare} и \ref{f:3poin}).
Определение инвариантов и их свойства составляют теорию, представленную в~\S\S\ref{s:wind}--\ref{s:minus2}.
Эта теория нужна также для теоремы~\ref{t:hocl-algor}.

\begin{example}\label{veplho-poin} Четыре замкнутые ломаные на рисунке~\ref{f:2points-a2-b3} попарно не гомотопны на плоскости без двух точек.
\end{example}

\begin{figure}
\centering
\compile{
\includegraphics[scale=0.3]{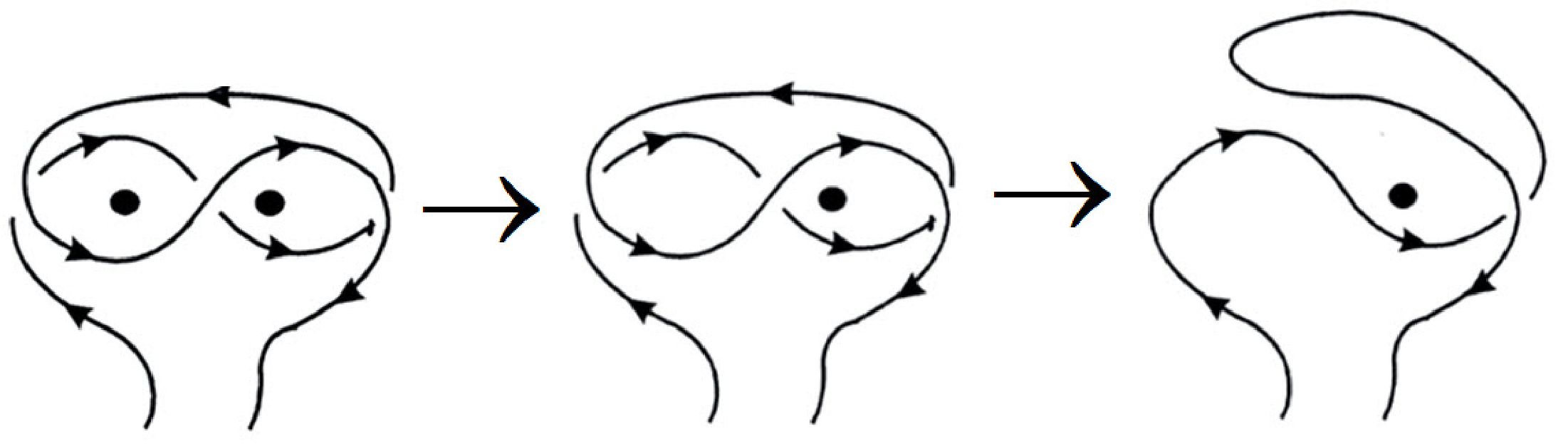}
}
\caption{После вынимания любого из гвоздей веревка падает}
\label{f:drop}
\end{figure}

\begin{example}[парадокс Пуанкаре]\label{veplho-borr} 
(a) На двух гвоздях, вбитых в~плоскую стену, можно так подвесить замкнутую веревку (с~тяжелой медалью), чтобы веревка не падала, но после вынимания любого гвоздя падала (см. левую часть рисунка~\ref{f:poincare} и рисунок~\ref{f:drop}).

Более строго, существуют две точки $P,Q$ и~замкнутая ломаная в~$\R^2-P-Q$, гомотопная одноточечной ломаной в~$\R^2-P$ и~в~$\R^2-Q$, но не в~$\R^2-P-Q$ (см. правую часть рисунка~\ref{f:poincare}).

(b) Аналогичное утверждение верно для трех гвоздей.
\end{example}

См. рисунки~\ref{f:poincare},~\ref{f:drop} и наглядную интерпретацию~в~\cite{Zi10},~\cite{DDM+}.
 
\begin{example}\label{p:3points}
Замкнутые ломаные на рисунке~\ref{f:3poin} не гомотопны на плоскости без трех точек.
\end{example}

Изучая настоящий текст, читатель освоит полезные и важные идеи и методы.
Они полезны и важны, поскольку дают следующий яркий результат, формулировка которого доступна неспециалисту в топологии (и многие другие результаты, см. учебники по топологии, например, \cite{Sk20}).

\begin{theorem}\label{t:hocl-algor} Для любого $n$ существует алгоритм распознавания гомотопности замкнутых ломаных на плоскости без $n$ точек.
\end{theorem}

Это конечная точка основных результатов по гомотопической классификации на плоскости без $n$ точек (мы начинаем с $n=1,2$, см. теоремы~\ref{pr:bij-1}, \ref{pr:bij-2} и~\ref{p:bij-n}). 
  
Вы увидите, как гомотопическая классификация связана с \emph{комбинаторикой слов}.
Естественно появится групповая структура, хотя и на несколько ином множестве \emph{отмеченных} замкнутых ломаных с точностью до \emph{отмеченной} гомотопности.
Это проливает свет на связь представленной классификации с \emph{комбинаторной теорией групп}.
Благодаря этой связи топологические методы могут быть использованы в теории групп.
Это отправная точка \emph{геометрической теории групп}~\cite{GG}.

\begin{figure}[h]\centering
\compile{
\includegraphics[scale=0.25]{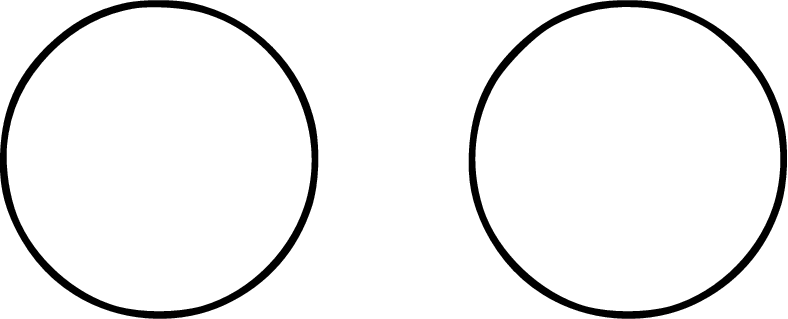} \qquad \includegraphics[scale=0.9]{3-1.eps}
}
\caption{Незацепленные и зацепленные кольца}\label{f:hopf}
\end{figure}

\begin{pr}[загадка]\label{3space}
В трехмерном пространстве даны два  

(a) незацепленных; \quad (b) зацепленных  

кольца, как показано на рисунке \ref{f:hopf}. 
Можно ли намотать и замкнуть веревку так, чтобы замкнутую веревку нельзя было оттащить далеко от двух колец, но можно было бы оттащить далеко от каждого кольца, разрезав другое кольцо? 

Строгая формулировка аналогична примеру~\ref{veplho-borr}. 
\end{pr}

\subsection*{О стиле этого текста}

\input{algor_style}

\section{Число оборотов: определение и обсуждение}\label{s:wind}

Всюду далее рассматриваемые точки и замкнутые ломаные расположены на плоскости.

Пусть $O,A,B$ "--- точки, причем $A\ne O$ и $B\ne O$ (но, возможно, $A=B$).
\emph{Ориентированным (или направленным) углом} $\angle AOB$ называется число  $t\in(-\pi,\pi]$, такое что вектор $\overrightarrow{OB}$ сонаправлен вектору, полученному из $\overrightarrow{OA}$ вращением на угол $t$ против часовой стрелки.
(Если Вы
можете рассматривать векторы на плоскости как комплексные числа, то можете переписать это условие как $\overrightarrow{OB}\upuparrows e^{it}\overrightarrow{OA}$.)
Далее рассматриваются ориентированные углы и слово <<ориентированный>> пропускается.

{\bf Числом оборотов} $w(l,O)$ замкнутой ломаной $l = A_1\ldots A_m$ вокруг не лежащей на ней точки $O$ называется количество оборотов при вращении вектора, начало которого находится в точке $O$, а конец обходит $l$ в заданном порядке вершин.
Строго говоря,
$$2\pi \cdot w(l,O) := \angle A_1OA_2+\angle A_2OA_3+\ldots+\angle A_{m-1}OA_m+\angle A_mOA_1.$$

\begin{figure}[ht]\centering
\compile{
\includegraphics[scale=0.65]{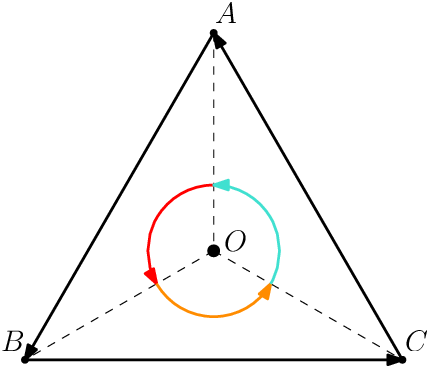}
\qquad
\includegraphics[scale=0.85]{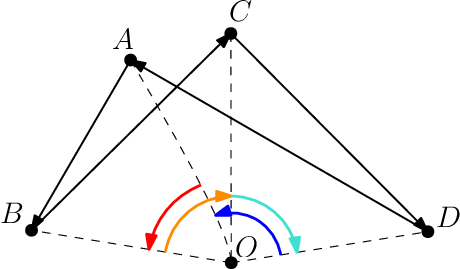}
}
\caption{$w(ABC,O) = +1 $ и $w(ABCD,O)=  0$}
\label{f:abco}
\end{figure}

Например, на рисунке~\ref{f:abco}
$w(ABC,O) = \dfrac{1}{2\pi} \left( \angle AOB + \angle BOC + \angle COA \right) = +1\quad\text{и}$
$$2\pi\cdot w(ABCD,O) =
\angle AOB + \angle BOC + \angle COD + \angle DOA =
\angle BOD + \angle DOB
= 0.$$

\begin{figure}[ht]\centering
\compile{
\includegraphics[scale=0.9]{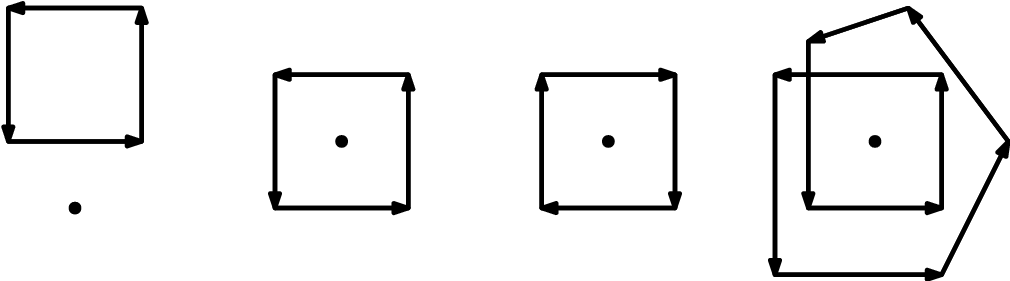}
}
\caption{Числа оборотов равны $0,~+1,~-1,~+2$}
\label{f:wn_examples}
\end{figure}

Пусть $ABC$ "--- треугольник и $O$ "--- точка внутри него. Тогда $w(ABCABC, O) = 2w(ABC) = \pm 2$.
Этот пример показывает, что числа оборотов для разных ломаных с одинаковым объединением их отрезков могут быть разными. 


\begin{proposition}\label{p:winsim} Число оборотов любого контура выпуклого многоугольника вокруг любой точки из внешности (внутренности) равно $0$ ($\pm 1$). 
См. рисунок~\ref{f:wn_examples}.
\end{proposition}

\begin{example}\label{p:surj-wn} Для любых целого числа $n$ и точки $O$ существует замкнутая ломаная, число оборотов которой вокруг $O$ равно $n$.
\end{example}

\begin{proposition}\label{p:noncl}
    Число оборотов $w(A_1\ldots A_m,O)$ является целым числом.
\end{proposition}

См. указания и доказательства в \cite[\S\S2,3]{AMS'}.
Подробнее о числе оборотов и связанных с ним понятиях см.~\cite{AMS'} и ссылки в этой работе.

\section{Классификация на плоскости без точки}
\label{s:minus}

В этом параграфе рассматриваются замкнутые ломаные и их гомотопность \emph{на плоскости без точки $O$}.

\begin{proposition}\label{p:winv}
Число оборотов вокруг~$O$ является инвариантом гомотопности. 
Иными словами, замкнутые ломаные с~разным числом оборотов вокруг~$O$ не го\-мо\-топны.
\end{proposition}

\textit{Набросок доказательства.} 
Достаточно доказать, что \textit{если замкнутая ломаная~$l_1$ получена из замкнутой ломаной~$l_2$ элементарным сокращением, то $w(l_1,O) = w(l_2,O)$}.

\begin{figure}[h]\centering
\compile{
\includegraphics[scale=0.3]{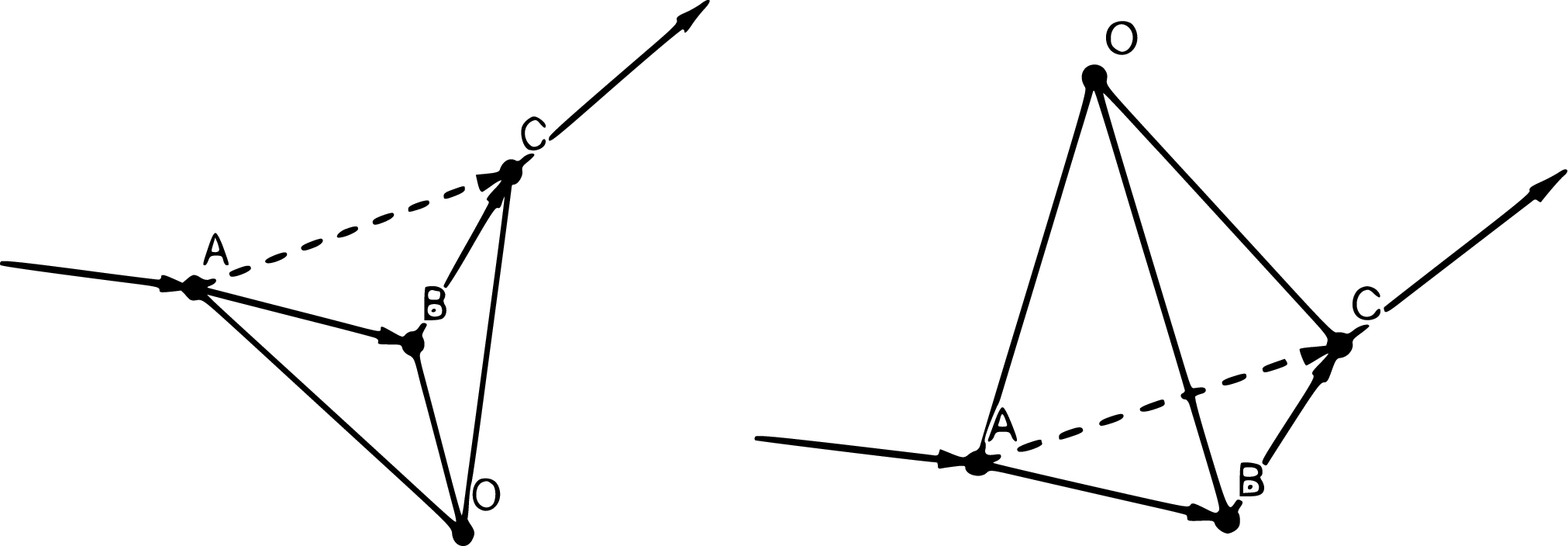}
}
\caption{Элементарное сокращение вершины $B$ и углы}
\label{f:wind}
\end{figure}

Докажем это.
Обозначим через $A,B,C$ последовательные точки ломаной $l_2$, такие что <<сокращение>> вершины $B$ дает ломаную~$l_1$.
Тогда
$$2\pi(w(l_1,O)-w(l_2,O)) = \angle AOB+\angle BOC-\angle AOC = 0,$$
где второе равенство справедливо, поскольку $O$ лежит вне треугольника~$ABC$ (см. рисунок~\ref{f:wind}; мы не доказываем соответствующее свойство ориентированных углов \cite[Утверждение 2.1]{AMS'}). 

\begin{proposition}\label{veplho-poly}
(a)~Замкнутая ломаная $l$ гомотопна одноточечной замкнутой ломаной, если  
$w(l,O)=0$.

(b)~Замкнутые ломаные с~одинаковым числом оборотов вокруг~$O$ го\-мо\-топны. 

Более точно, возьмем контур $\partial\Delta$ треугольника $\Delta$, содержащего внутри себя точку~$O$.
Тогда любая замкнутая ломаная $l$ гомотопна $|w(l,O)|$"=кратному обходу контура $\partial\Delta$

$\bullet$ против часовой стрелки при $w(l,O) \geqslant 0$,

$\bullet$ по часовой стрелке при $w(l,O) \leqslant 0$.
\end{proposition}

Наброски доказательств этого и других утверждений данного параграфа приведены в его конце. 

Следующая теорема~\ref{pr:bij-1} следует из примера~\ref{p:surj-wn} и утверждений~\ref{p:winv}, \ref{veplho-poly}.b.

\begin{theorem}\label{pr:bij-1} 
Следующее отображение является взаимно однозначным соответствием между множеством гомотопических классов замкнутых ломаных на плоскости без точки и множеством целых чисел. 
Отображение переводит замкнутую ломаную (точнее, ее гомотопический класс) в число ее оборотов вокруг точки.
\end{theorem}

Из утверждения \ref{veplho-poly}.b видно, что описанная выше двумерная теорема~\ref{pr:bij-1} по сути является одномерной.
Аналогичное понимание используется для гомотопической классификации на плоскости без более чем одной точки.
 
Луч на плоскости называется \emph{лучом общего положения относительно замкнутой ломаной}, если он не проходит через ее вершины.

\begin{example}
Существует замкнутая ломаная $l$ и луч общего положения $OP$ относительно $l$, такие что количество точек пересечения ломаной $l$ и луча $OP$ имеет 

(a) ту же четность, что и $w(l, O)$; \quad (b) другую четность, чем $w(l, O)$.
\end{example}

\begin{proposition}\label{p:ray}
Возьмем замкнутую ломаную $l$ и луч общего положения $OP$ относительно $l$. 
Количество звеньев ломаной $l$, пересекающих луч $OP$, имеет ту же четность, что и $w(l, O)$.
\end{proposition}

Для доказательства утверждения~\ref{p:ray}, использующего гомотопию, полезны утверждения~\ref{veplho-poly}.b и \ref{p:apprgen}. 
Доказательство утверждения~\ref{p:ray}, не использующее гомотопию, см. в~\cite[\S4]{AMS'}. 

\begin{proposition}\label{p:apprgen}
Пусть даны две гомотопные замкнутые ломаные. 
Возьмем луч общего положения $OP$ относительно каждой из этих ломаных.
Тогда эти ломаные могут быть соединены гомотопией, относительно каждой из ломаных которой луч $OP$ является лучом общего положения. 
\end{proposition}


\begin{figure}[h]\centering
\compile{
\includegraphics[scale=0.25]{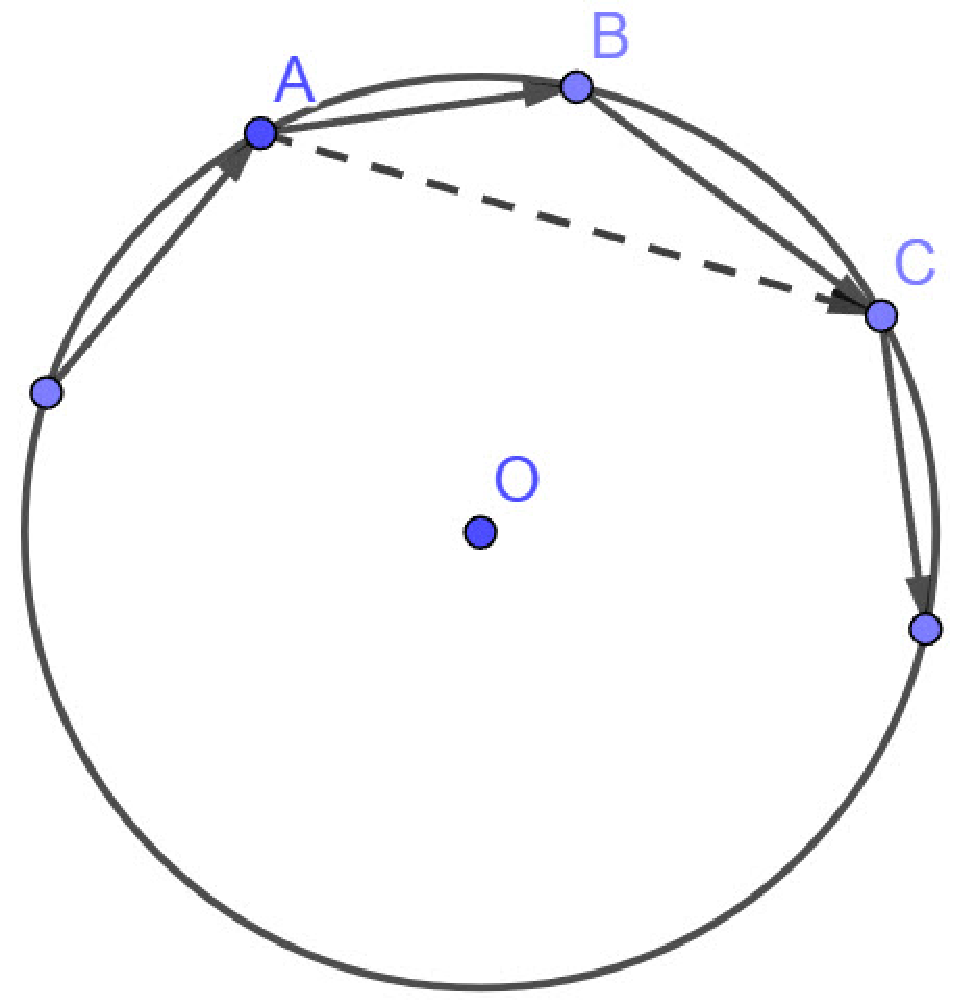}
}
\caption{Элементарное сокращение обхода}
\label{f:rou}
\end{figure}

Приведем в виде задачи версию теоремы~\ref{pr:bij-1} на более традиционном языке.

\begin{pr}\label{veplho-polyc}
Назовем \emph{обходом} ориентированную циклическую последовательность точек (фиксированной) окружности, в~которой никакие две последовательные точки не являются диаметрально противоположными.
\emph{Элементарным сокращением} обхода называется удаление такой его точки~$B$, что для соседних с~$B$ точек $A,C$ выпуклая оболочка точек $A,B,C$ не содержит центра $O$ окружности (см. рисунок~\ref{f:rou}).
Два обхода называются \emph{гомотопными}, если они могут быть соединены последовательностью обходов, в~которой один из любых двух последовательных обходов получен из другого элементарным сокращением.

(a)~Обход~$l$ гомотопен одноточечному обходу тогда и~только тогда, когда $w(l,O)=0$.

(b)~(загадка)
Опишите обходы с~точностью до гомотопности.
\end{pr}



\begin{proposition}\label{p:3rays} 
    Пусть $A_0A_1A_2$ "--- правильный треугольник, а точка $O$ "--- его центр.
    Для $j=0,1,2$ обозначим через $l_j$ ломаную, не пересекающую луч $OA_j$ и соединяющую $A_{j+1}$ с $A_{j+2}$, где нумерация берется по модулю 3.
    Тогда $w(l_0l_1l_2, O) = \pm 1$.  
\end{proposition}

\subsection*{Наброски доказательств некоторых утверждений}~

\medskip
\noindent\textbf{\ref{veplho-poly}.a.} 
Достаточно доказать, что \textit{если замкнутая ломаная имеет более одной вершины и~имеет нуль оборотов, то к~ней можно применить элементарное сокращение}.

Докажем это.
Если в~замкнутой ломаной некоторые две последовательные вершины совпадают, то можно одну из них элементарно сократить.
Пусть теперь любые две последовательные вершины различны.
Так как сумма углов равна нулю, то найдутся два соседних угла $\angle AOB$ и~$\angle BOC$ разных знаков.
Тогда точки~$A$ и~$C$ лежат в~одной полуплоскости относительно прямой~$BO$.
Значит, треугольник~$ABC$ не содержит точку~$O$.
Поэтому можно сделать элементарное сокращение вершины~$B$.

\begin{figure}[h]\centering
\compile{
\includegraphics[scale=.3]{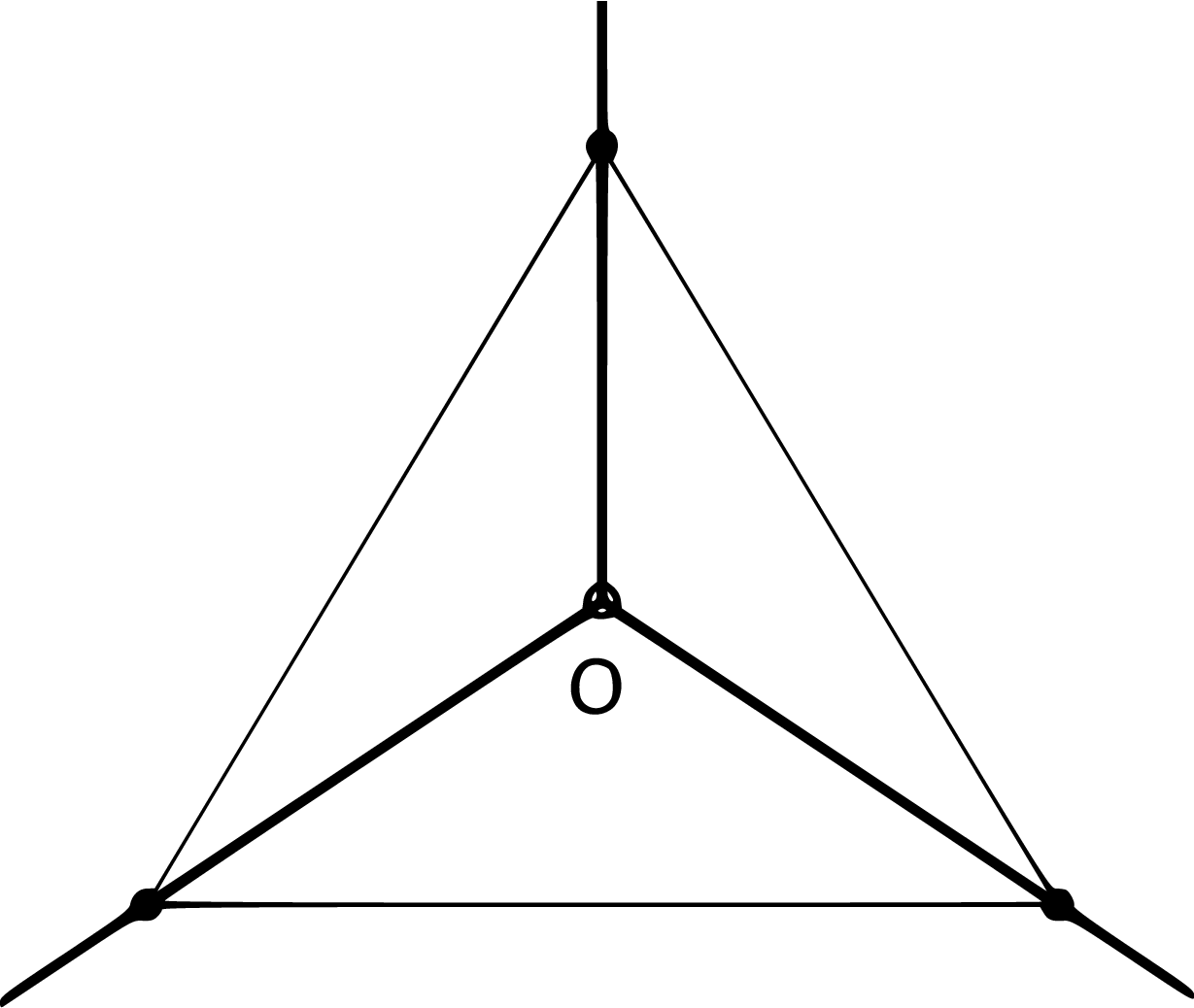}
}
\caption{Дополнительные лучи и <<двойственный>> треугольник на плоскости без точки}
\label{f:compl1}
\end{figure}

\medskip
\noindent\textbf{\ref{veplho-poly}.b.}
Достаточно доказать, что \textit{любая замкнутая ломаная гомотопна замкнутой ломаной, содержащей только вершины треугольника $\Delta$} (рисунок~\ref{f:compl1}).

Докажем это.
Назовем пару $\{ XY,\ B \}$, состоящую из отрезка $XY$ ломаной $l$ и вершины $B$ треугольника $\Delta$, \emph{неудобной} (для $l$), если $B \notin \{X, Y\}$ и $B \in \angle XOY$.
Назовем ломаную \emph{удобной}, если 

$\bullet$ она содержит по крайней мере одну вершину треугольника $\Delta$, и 

$\bullet$ ни для каких звена $XY$ ломаной $l$ и вершины $B$ треугольника $\Delta$ пара $\{XY, B \}$ не является неудобной.
   
В следующем абзаце мы доказываем, что \emph{любая замкнутая ломаная $l$ гомотопна удобной.} 

Можно считать, что $l$ не является удобной. 
Сделав (если необходимо) несколько обратных операций к элементарному сокращению, мы можем считать, что $l$ содержит вершину треугольника $\Delta$.
Тогда существует неудобная пара $\{ XY,\ B \}$ для $l$.
Добавим $B$ к $l$ между $X$ и $Y$, используя обратную операцию к элементарному сокращению. 
Количество неудобных пар для полученной замкнутой ломаной меньше, чем для $l$.
С помощью конечного числа таких операций получим удобную замкнутую ломаную, гомотопную $l$.

Любую вершину удобной замкнутой ломаной, отличную от вершин треугольника $\Delta$, можно удалить элементарным сокращением.
Полученная замкнутая ломаная также удобна.
Следовательно, с помощью конечного числа таких элементарных сокращений мы получим замкнутую ломаную, содержащую только вершины треугольника $\Delta$.

\medskip
\noindent\emph{Набросок другого доказательства утверждения}~\textbf{\ref{veplho-poly}.b} \emph{(А. Абзалилов).}
Выберем вершину $A$ ломаной $l$.  
Обозначим через $\overline l$ (не циклическую) последовательность точек, начинающуюся с $A$ соответствующую ломаной $l$.
Обозначим через $\overline d=B_1B_2B_3$ (не циклическую) последовательность вершин треугольника $\Delta$,   ориентированную так, что числа $w(\overline d,O)$ и $w(l,O)$ имеют одинаковый знак.
Можно считать, что $O$ не лежит на прямой $AB_1$.
Обозначим  
$$n:=|w(l,O)|, \quad \overline d^{-1}:=B_3B_2B_1, \quad \overline d^{-n}:= \underbrace{\overline d^{-1} \ldots \overline d^{-1}}_{n}\quad\text{и}\quad x:=\overline l A B_1 \overline d^{-n} A.$$
Аналогично определим $\overline d^n$.
Теперь нужная гомотопия от $\overline l$  к $n$-кратному обходу $\overline d^{n}$ контура $\partial\Delta$ строится так:
$$\overline l \ \rightarrow\ \overline lA \ \rightarrow\ \overline lA B_1\ \rightarrow\ \overline lAB_1B_1 \ 
\rightarrow \ \overline lAB_1B_3B_3B_1\ \rightarrow$$ 
$$\overline lAB_1B_3B_2B_2B_3B_1 \rightarrow\ \overline lAB_1B_3B_2B_1B_1B_2B_3B_1 \ \rightarrow$$ 
$$\ \ldots\ \rightarrow\ \overline lAB_1\overline d^{-n}\overline d^{n} B_1\  \rightarrow\ \overline lAB_1\overline d^{-n}A\overline d^{n} B_1\ =
\ x\overline d^{n}B_1\ \Rightarrow\ A\overline d^{n}B_1\ \rightarrow\ \overline d^{n} B_1\ \rightarrow
\ \overline d^{n}.$$
Здесь все элементарные сокращения `$\rightarrow$' возможны (так как точка $O$ не лежит на отрезке $A B_1$ 
и $O\not\in\partial\Delta$), а гомотопия `$\Rightarrow$' строится так.  
Поскольку $w(x,O)=0$, ломаная $x$ гомотопна одноточечной ломаной.
Более того, существует гомотопия между $x$ и одноточечной ломаной $A$
у каждой промежуточной ломаной которой есть вершина $A$.
Припишем к каждой ломаной в этой гомотопии последовательность $\overline d^{n} B_1$.
Получим гомотопию `$\Rightarrow$'.

\medskip
\noindent\textbf{\ref{p:apprgen}.}
(Этот набросок получен редактированием текста А. Мизева.)

Возьмем любую гомотопию между заданными ломаными $l_1,l_2$.
Если ни одна из вершин ломаных гомотопии не лежит на луче $U := OP$,
то утверждение доказано. 
Иначе обозначим через $A \in U$ такую вершину. 
Точку $A'\not\in U$ назовем \emph{похожей на $A$}, если для любого <<треугольника сокращения>> $ABC$ треугольник $A'BC$ не содержит точку $O$.

В следующем абзаце докажем, что существует точка $A\not\in U$, похожая на $A$.  

Используем следующее утверждение. 
\emph{Пусть треугольник $XYZ$ не содержит точку $O$. 
Тогда существует $\varepsilon>0$ такое, что для любой точки $X'$, удаленной от $X$ не более, чем на $\varepsilon$, треугольник $X'YZ$ не содержит $O$.}
Применим его к каждому <<треугольнику сокращения>>, одной из вершин которого является $A$.
Возьмем $\varepsilon_1>0$, меньшее всех полученных $\varepsilon$'в. 
Существует точка, удаленная от $A$ не более чем на $\varepsilon_1$, и не лежащая  в $U$. 
Она похожа на $A$. 

В каждой ломаной гомотопии, содержащей $A$, заменим каждую вершину $A$ на похожую вершину $A'$. 
Получим гомотопию между $l_1$ и $l_2$.
Так мы уменьшим количество вершин в ломаных гомотопии, лежащих в $U$. 
Будем повторять это, пока существует вершина одной из ломаных гомотопии, лежащая в $U$. 
В итоге получим нужную гомотопию между $l_1$ и $l_2$.

\medskip
\noindent\textbf{\ref{p:3rays}.}
Для каждого $j = 0, 1, 2$ каждую вершину ломаной $l_j$ спроецируем гомотопией вне $O$ на прямую $A_{j+1}A_{j+2}$.
Получим ломаную $m_j$ и гомотопию вне $O$ ломаной $l_0l_1l_2$ к ломаной $m_0m_1m_2$.
По утверждению~\ref{p:winv} имеем $w(l_0l_1l_2, O) = w(m_0m_1m_2, O)$.
Элементарно сократим все вершины ломаной $m_0m_1m_2$, кроме вершин $A_0, A_1, A_2$.
Получим гомотопию вне $O$ ломаной $m_0m_1m_2$ к треугольнику $A_0A_1A_2$.
По утверждению~\ref{p:winv} имеем $w(m_0m_1m_2, O) = w(A_0A_1A_2, O) = \pm 1$.

\medskip
См. другое доказательство утверждения~\ref{p:3rays}, не используещее понятие гомотопии, в~\cite[конец \S3]{AMS'}.

\section{<<Разрешение>> парадокса Пуанкаре}\label{s:poipar}

Очевидно, что в примере~\ref{veplho-borr}.a, приведенном на рисунке~\ref{f:poincare}, числа оборотов  вокруг точек $P$ и $Q$ равны нулю.
Негомотопность в примере~\ref{veplho-borr}.a следует из утверждения~\ref{p:invar-link} (или утверждения~\ref{p:invar-2eco}). 
Идея заключается в использовании одного из следующих инвариантов гомотопности.

\begin{figure}[h]\centering
\compile{
\includegraphics[scale=0.9]{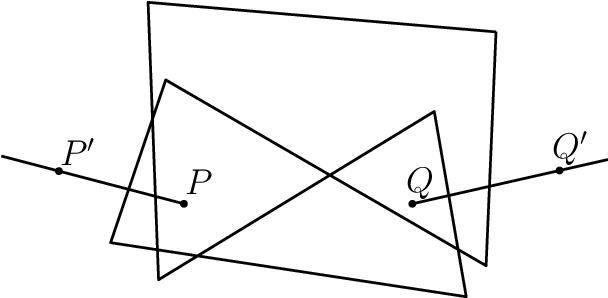}
\quad
\includegraphics[scale=0.4]{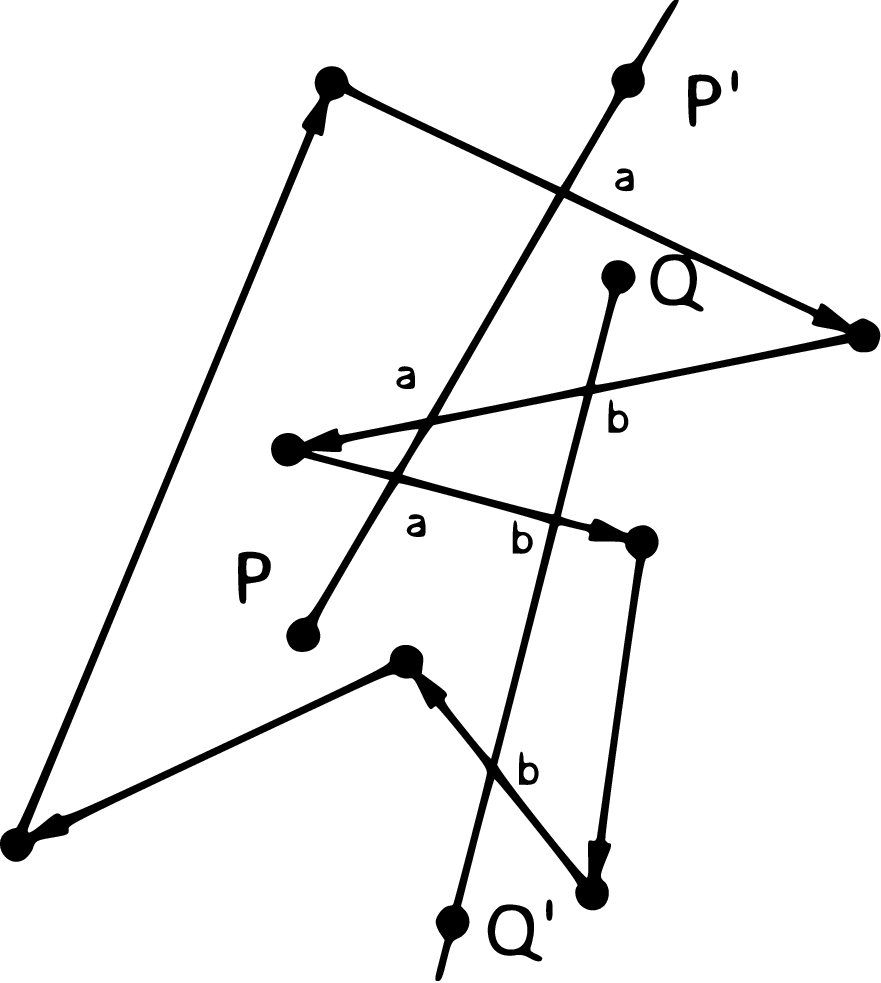}
}
\caption{Циклические слова Пуанкаре по модулю 2: $abab$ (слева), $abaabb$ (справа)}
\label{f:cycl}
\end{figure}

В \S\S\ref{s:poipar}--\ref{s:minus2}

$\bullet$ $l$ "--- замкнутая ломаная на плоскости без двух точек $P,Q$;

$\bullet$ $PP'$ и $QQ'$ "--- непересекающиеся лучи общего положения относительно $l$.

Будем двигаться по $l$ и выписывать букву $a$ (букву $b$) при каждом пересечении луча $PP'$ (луча $QQ'$).
Полученное ориентированное циклическое слово называется \emph{циклическим словом Пуанкаре по модулю 2} 
(см. рисунок~\ref{f:cycl}). 

\begin{example}
\label{p:poin2-a}
Любое циклическое слово из букв $a$ и $b$ является циклическим словом Пуанкаре по модулю 2 некоторой замкнутой ломаной на плоскости без двух точек (для некоторых лучей $PP'$, $QQ'$).
\end{example}

\begin{proposition}
\label{p:poin2-b}
Четность числа оборотов замкнутой ломаной $l$ вокруг точки $P$ (точки $Q$) равна четности количества букв $a$ (букв $b$) в циклическом слове Пуанкаре по модулю 2 ломаной $l$ (для любых лучей $PP'$, $QQ'$).
\end{proposition}

\begin{figure}[h]\centering
\compile{
\includegraphics[scale=1.]{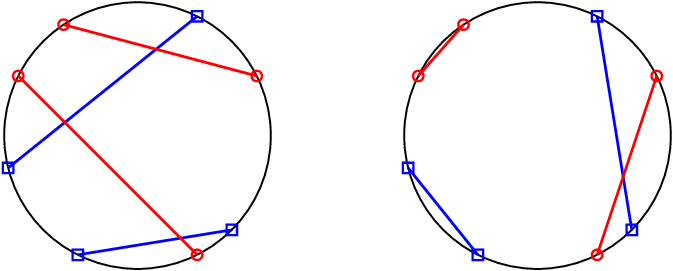}
}
\caption{Два разбиения на пары четырех красных точек (маленькие кружочки) и четырех синих точек (маленькие квадраты) на окружности}
\label{f:circles}
\end{figure}

Пусть на окружности дано четное число красных точек и четное число синих точек, причем все эти точки попарно различны.
Разобьем красные точки на пары, и синие точки на пары.
Соединим точки в каждой красной (синей) паре красной (синей) хордой (см. рисунок~\ref{f:circles}).
Множества красных точек и синих точек называются \emph{зацепленными}, если количество пар пересекающихся красной и синей хорд нечетно (ср.~\cite[утверждение 4.9.5a]{Sk}).

\begin{proposition}
\label{p:mil-dis} 
Зацепленность не зависит от выбора разбиений на пары (см. рисунок~\ref{f:circles}).
\end{proposition}

Назовем \emph{интересной} замкнутую ломаную $l$ в плоскости без точек $P, Q$, для которой $w(l, P)$ и $w(l, Q)$ четны. 
По утверждению~\ref{p:winv} интересность сохраняется при замене ломаной на гомотопную.  
По утверждению~\ref{p:poin2-b} количество букв $a$ (букв $b$) в циклическом слове Пуанкаре по модулю 2 интересной замкнутой ломаной четно.
Назовем интересную замкнутую ломаную \emph{зацепленной}, если буквы $a$ и $b$ в ее циклическом слове Пуанкаре по модулю 2 зацеплены. 

\begin{proposition}\label{p:invar-link} 
Зацепленность интересной замкнутой ломаной не меняется при

(a) изменении лучей $PP'$ и $QQ'$; 

(b) замене замкнутой ломаной $l$ на гомотопную. 
\end{proposition}

\begin{proof}[Набросок вывода п.~(b) из п.~(a)]
(Этот набросок получен редактированием текста А. Абзалилова.)
Достаточно доказать, что зацепленность по модулю~2 не меняется при элементарном сокращении.
По п.~(a) можно взять лучи $PP'$ и $QQ'$ общего положения с ломаной до элементарного сокращения, так чтобы они не пересекали друг друга, и ни один из них не пересекал <<треугольника сокращения>>.
Тогда элементарное сокращение не меняет циклическое слово Пуанкаре по модулю~2.
Значит оно не меняет и зацепленности по модулю~2.  
\end{proof}

В п.~(b) циклическое слово Пуанкаре по модулю 2 может не быть определено для некоторых промежуточных ломаных гомотопии.
Поэтому аналог утверждения~\ref{p:apprgen} для плоскости без двух точек полезен для прямого доказательства утверждения~\ref{p:invar-link}.b.

\section{Инвариант для плоскости без двух точек}\label{s:poiword}

\emph{Данный параграф подводит к \S\ref{s:minus2}. Тем не менее, этот параграф можно пропустить, поскольку он не используется в дальнейшем.}

\begin{example}\label{p:link-noncomp} Существует ломаная $l$, такая что $w(l, P) = w(l, Q) = 0$, ломаная $l$ не зацеплена, но $l$ не гомотопна одноточечной замкнутой ломаной. 
\end{example}

Рассмотрим множество всех (конечных) циклических слов (включая пустое слово) из букв $a,b$.
Для такого слова \emph{элементарное сокращение} "--- замена любого из подслов $aa,bb$ пустым подсловом.
Такое слово называется \emph{экономным}, если в нем нет ни одного из подслов $aa$ и $bb$.
 
\begin{proposition}\label{p:wd-2eco} 
Циклическое слово из букв $a,b$ посредством элементарных сокращений дает единственное экономное слово.
\end{proposition}

\begin{proof}[Набросок доказательства (Е. Волович)]
Пусть, напротив, существует циклическое слово $A$, из которого элементарными сокращениями можно получить два разных экономных слова $E_1$ и $E_2$.

Фразу <<слово $Y$ получено некоторым элементарным сокращением слова $X$>> сокращаем до $X\to Y$, а фразу <<слово $Y$ получено несколькими элементарными сокращениями слова $X$>> сокращаем до $X\Rightarrow Y$. 

Пока существует слово $B$, из которого элементарными сокращениями можно получить два разных экономных слова, и  $A\Rightarrow B$, будем заменять $A$ на $B$.
Этот процесс конечен, так как с каждой заменой длина слова $A$ уменьшается.

После этого для $i=1,2$ обозначим через $B_i$ слово, такое что $A\to B_i$ и $B_i\Rightarrow E_i$. 
Тогда $B_1\ne B_2$.  
Поэтому пары букв в $A$, сокращением которых получаются $B_1$ и $B_2$, не пересекаются.
Значит, если сократить обе пары последовательно, то получится такое слово $C$, что $B_1\to C$ и $B_2\to C$.

Так как $A\Rightarrow C$, то из $C$ элементарными сокращениями можно получить только одно экономное слово $E$. 
Так как $A\to B_1$, $B_1\Rightarrow E_1$ и $B_1\Rightarrow E$, то $E=E_1$. 
Аналогично $E=E_2$.
Противоречие.
\end{proof}

\emph{Экономная форма} (или нормальная форма) циклического слова $w$ из букв $a, b$ — это экономное слово, полученное из $w$ элементарными сокращениями.
Обозначим через $E_{2,2,c}$ множество всех экономных слов.
Обозначим через $e_2(l)\in E_{2,2,c}$ экономную форму циклического слова Пуанкаре по модулю 2 ломаной $l$. 
Априори $e_2(l)$ зависит от лучей $PP'$ и $QQ'$.

\begin{proposition}[ср. с утверждением~\ref{p:invar-link}]\label{p:invar-2eco} Слово $e_2(l)$ не изменяется при 

(a) изменении лучей $PP'$ и $QQ'$;  \quad 

(b) замене замкнутой ломаной $l$ на гомотопную.  
\end{proposition}

\begin{example}\label{p:ncomp-2eco}  
Существует ломаная $l$, не гомотопная одноточечной замкнутой ломаной, но для которой слово $e_2(l)$ пусто.
\end{example}

\begin{proof}[Набросок построения (получен редактированием текста А. Абзалилова)]
Отметим на плоскости точки $O, P_1, P_2$ так, чтобы $O$ была серединой отрезка $PQ$, никакие из отмеченных точек не лежали на одной прямой, и точка $P$ лежала внутри треугольника $OP_1P_2$. 
Контрпримером является замкнутая ломаная $l:=(OP_1P_2)^2$.

Модуль количества ее оборотов вокруг точки $P$ равен 2.
  Значит, ломаная $l$ не гомотопна одноточечной.

  Направим луч $QQ'$ так, чтобы он не пересекал контура треугольника $OP_1P_2$.
  Тогда циклическое слово Пуанкаре по модулю 2 равно $aa$.
  Поэтому слово $e_2(l)$ пусто.
\end{proof}

Два циклических слова из букв $a,b$ называются \textit{эквивалентными}, если их можно соединить последовательностью циклических слов, в которой одно из любых двух последовательных циклических слов может быть получено из другого элементарным сокращением.
Обозначим через $F_{2,2,c}$ множество всех классов эквивалентности.

\begin{proposition}\label{p:word2-example}
(a) Слово $abab$ не эквивалентно пустому слову.

(b) Существуют ли два неэквивалентных циклических слова, оба состоящие из нечетного числа букв $a$ и четного числа букв $b$?

(c) Никакие два различных экономных слова не эквивалентны.

(Это следует из утверждения~\ref{p:wd-2eco}.)

(d) (загадка) Постройте <<естественное>> взаимно однозначное соответствие  $E_{2,2,c}\to F_{2,2,c}$.

(Указание: используйте п.~(c).)

(e) (загадка) Опишите $E_{2,2,c}$ или, что то же самое, $F_{2,2,c}$.
\end{proposition}

\begin{proposition}[ср. с утверждением~\ref{p:invar-2eco}]
\label{p:invar-2eq}
Класс эквивалентности в $F_{2,2,c}$ циклического слова Пуанкаре по модулю 2 не меняется при

(a) изменении лучей $PP'$ и $QQ'$;  \quad 

(b) замене замкнутой ломаной $l$ на гомотопную.  
\end{proposition}

\section{Классификация на плоскости без двух точек}
\label{s:minus2}

Приведенные выше гомотопические инварианты не являются полными (см. примеры~\ref{veplho-borr},~\ref{p:3points}, \ref{p:link-noncomp},~\ref{p:ncomp-2eco}).
Теперь мы готовы ввести полный инвариант.

В этом параграфе рассматриваются замкнутые ломаные и их гомотопия \emph{на плоскости без двух точек $P$ и $Q$}.

\begin{figure}[h]\centering
\compile{
\includegraphics[scale=.8]{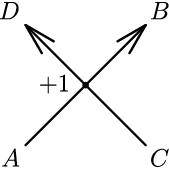}\quad 
\includegraphics[scale=.8]{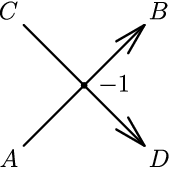}
}
\caption{Знак точки пересечения}
\label{f:sign}
\end{figure}

Пусть $A,B,C,D$ "--- точки на плоскости, никакие три из которых не лежат на одной прямой.
{\bf Знаком} точки пересечения ориентированных отрезков $\overrightarrow{AB}$ и $\overrightarrow{CD}$ на плоскости
назовем $+1$, если обход $ABC$ происходит по  часовой стрелке, и $-1$ в противном случае (см. рисунок~\ref{f:sign}).

\begin{proposition}[ср. с утверждением~\ref{p:ray}]\label{p:sign-ray}
Возьмем замкнутую ломаную $m$ на плоскости без точки $O$ и луч $OR$ общего положения относительно $m$. 
Тогда $w(m,O)$ равно сумме по звеньям $x$ ломаной $m$ знаков точек пересечения звена $x$ и луча $OR$.
\end{proposition}

Аналогично~\S\ref{s:poipar} определим для замкнутой ломаной и пары лучей $PP'$ и $QQ'$ общего положения \emph{циклическое слово Пуанкаре} из букв $a,b,a^{-1},b^{-1}$.

\begin{figure}[h]\centering
\compile{
\includegraphics[scale=0.5]{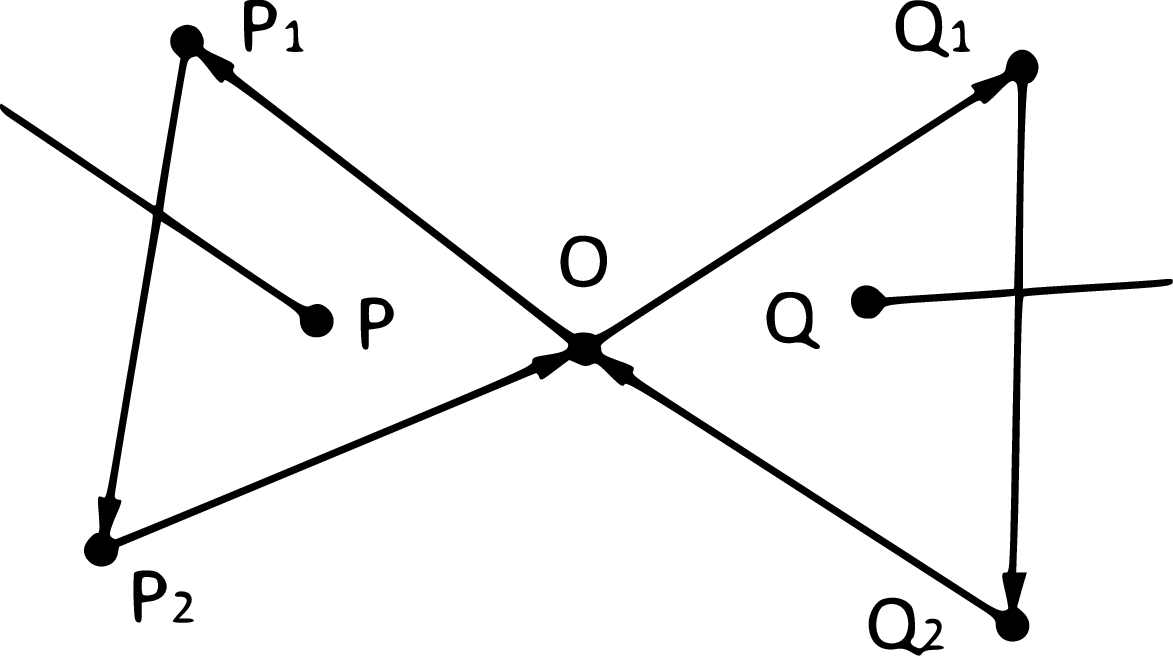}
}
\caption{<<Восьмерка>> и лучи общего положения}
\label{f:eight}
\end{figure}

\begin{example}[ср. с примером~\ref{p:poin2-a}]
\label{p:cycpoi}
Любое циклическое слово из букв $a,b,a^{-1},b^{-1}$ является циклическим словом Пуанкаре некоторой замкнутой ломаной
(для некоторых лучей $PP'$, $QQ'$; см. рисунок \ref{f:eight}).
\end{example}

Рассмотрим множество всех (конечных) циклических слов (включая пустое слово) из букв $a,a^{-1},b,b^{-1}$. 
Для такого слова \emph{элементарное сокращение} "--- замена любого из подслов $aa^{-1},a^{-1}a,bb^{-1},b^{-1}b$ пустым подсловом.
Такое слово называется \emph{экономным}, если у него нет ни одного из этих подслов.

\begin{proposition}
\label{p:wd-eco} 
Циклическое слово из букв $a,b, a^{-1}, b^{-1}$ посредством элементарных сокращений дает единственное экономное слово.
\end{proposition}

\textit{Экономная форма} (или нормальная форма) циклического слова $w$ из букв $a,b, a^{-1}, b^{-1}$ — это экономное слово, полученное из $w$ элементарными сокращениями. 
Обозначим через $E_{2,c}$ множество всех циклических экономных слов. 
Обозначим через $e(l)\in E_{2,c}$ экономную форму циклического слова Пуанкаре ломаной $l$. 
Априори $e(l)$ зависит от лучей $PP'$ и $QQ'$.

\begin{proposition}[ср. с  утверждением~\ref{p:invar-2eco}]\label{p:invar-eco}
Слово $e(l)$ не изменяется при

(a) изменении лучей $PP'$ и $QQ'$;  \quad 

(b) замене замкнутой ломаной $l$ на гомотопную.  
\end{proposition}

Для доказательства п. (b) полезен аналог утверждения~\ref{p:apprgen} для плоскости без двух точек.

\begin{proposition}[ср. с примером~\ref{p:ncomp-2eco}]\label{p:comp-eco}
Если $e(l_1)=e(l_2)$ для замкнутых ломаных $l_1$ и $l_2$, то $l_1$ и $l_2$ гомотопны. 
\end{proposition}

\begin{figure}[h]\centering
\compile{
\includegraphics[scale=.3]{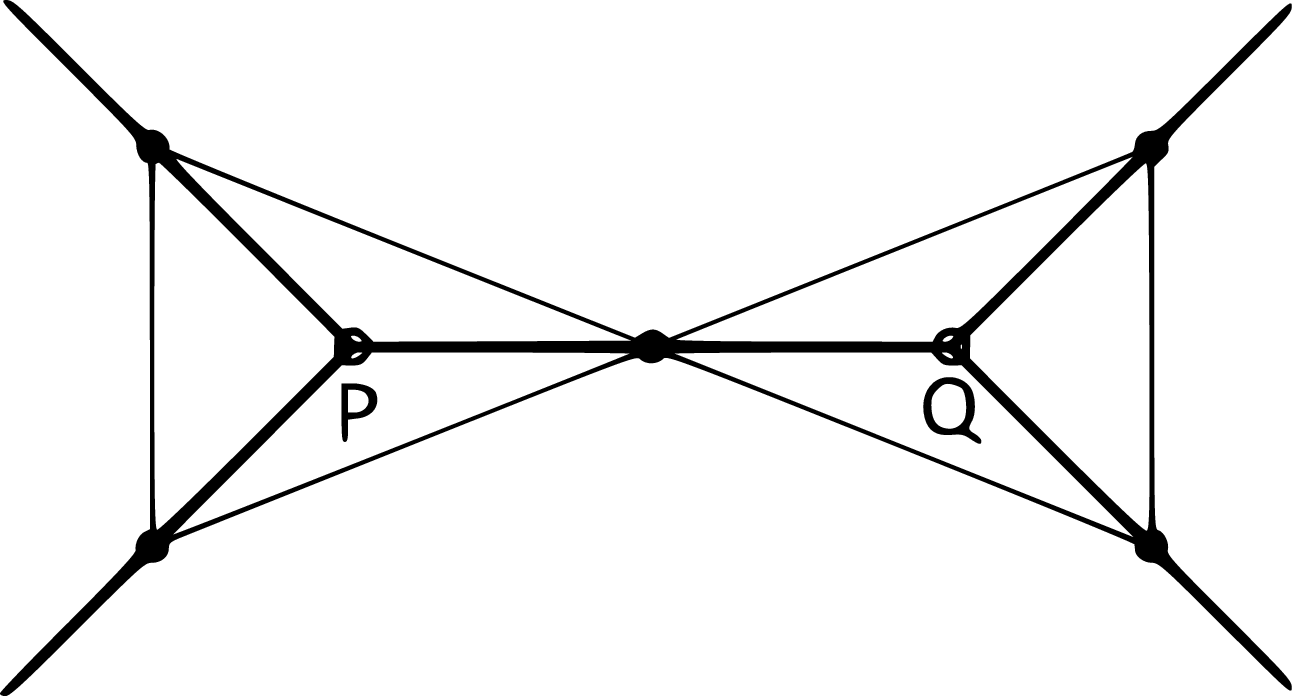}
}
\caption{Дополнительное построение и <<двойственная>> восьмерка на плоскости без двух точек $P,Q$}
\label{f:compl2}
\end{figure}

См. рисунок~\ref{f:compl2} и набросок доказательства теоремы~\ref{p:bij-n} в конце параграфа. 
Ср. задачу~\ref{p:defret}. 

Следующая теорема~\ref{pr:bij-2} следует из примера~\ref{p:cycpoi} и утверждений~\ref{p:invar-eco},~\ref{p:comp-eco}.

\begin{theorem}\label{pr:bij-2} 
Следующее отображение является взаимно однозначным соответствием между множеством гомотопических классов замкнутых ломаных на плоскости без двух точек и множеством $E_{2,c}$.
Отображение переводит замкнутую ломаную в экономную форму ее циклического слова Пуанкаре.
\end{theorem}

Из этой теоремы вытекает теорема~\ref{t:hocl-algor} для $n = 2$.

\begin{remark}\label{r:equiv}
Аналогично \S\ref{s:poiword}, определим отношение эквивалентности на множестве циклических слов из букв $a,b,a^{-1},b^{-1}$. 
Обозначим через $F_{2,c}$ множество всех классов эквивалентности.
Аналогично утверждению~\ref{p:word2-example}, существует <<естественное>> взаимно однозначное соответствие $E_{2,c}\to F_{2,c}$.
\end{remark}

\begin{pr}[загадка]\label{p:defret} (a) Определите гомотопность ориентированных циклов в графе <<восьмерка>>.
Постройте взаимно однозначное соответствие между гомотопическими классами замкнутых ломаных на плоскости без двух точек и гомотопическими классами ориентированных циклов в графе <<восьмерка>>. 

(b) Верен аналог теоремы~\ref{pr:bij-2} для графа <<восьмерка>> вместо плоскости без двух точек.
\end{pr}

\section{Классификация на плоскости без $n$ точек}
\label{s:minusn}

Обозначим выколотые точки через $P_1, P_2, \ldots, P_n$. 
Произвольный набор из $n$ попарно непересекающихся лучей с началами в выколотых точках назовем \emph{разрезом}.
Разрез называется \emph{разрезом общего положения относительно замкнутой ломаной}, если каждый луч разреза является лучом общего положения относительно этой ломаной.

Возьмем разрез $r_1, r_2, \ldots, r_n$ общего положения относительно замкнутой ломаной.
Будем двигаться по ломаной и выписывать букву $p_i$ (букву $p_i^{-1}$) при каждом пересечении звена ломаной с лучом $r_i$, если знак точки их пересечения равен $+1$ ($-1$).
Полученное ориентированное циклическое слово  называется \emph{циклическим словом Пуанкаре 
(для разреза $r_1, \ldots, r_n$)}.

\emph{Экономные циклические слова} и \emph{экономная форма} циклического слова из букв $p_1,\ldots,p_n$, $p_1^{-1},\ldots,p_n^{-1}$ определяются аналогично случаю $n=2$ (\S\ref{s:minus2}).

\begin{example}[см. левую часть рисунка~\ref{f:3poin}]\label{e:3raysn}
Для некоторой замкнутой ломаной $l$ на плоскости без трех точек (без $n \geqslant 3$ точек) существуют два разреза, таких что экономные формы циклических слов Пуанкаре ломаной $l$ для этих разрезов различны.
\end{example}


Этот пример мотивирует следующее определение.
Зафиксируем набор из $n$ попарно непересекающихся углов $\alpha_1, \ldots, \alpha_n$ с вершинами $P_1, \ldots, P_n$, соответственно.
Обозначим через $e(l) := e(l; \alpha_1, \ldots, \alpha_n)$ экономную форму циклического слова Пуанкаре ломаной $l$ для разреза $r_1, \ldots, r_n$ общего положения относительно $l$, лучи которого лежат в $\alpha_1, \ldots, \alpha_n$, соответственно.
Априори $e(l)$ зависит от разреза $r_1, \ldots, r_n$.


\begin{proposition}[ср. с утверждением \ref{p:invar-eco}]
\label{p:invar-eco-n}
Слово $e(l)$ не изменяется при

(a) изменении разреза $r_1, \ldots, r_n$;  \quad 

(b) замене замкнутой ломаной $l$ на гомотопную.  
\end{proposition}

\noindent\emph{Указание к п. (a).} Достаточно доказать утверждение для двух разрезов, отличающихся только в одном индексе $i = 1, \ldots, n$.
Дальнейшее рассуждение аналогично доказательству утверждения \ref{p:invar-eco}.a, поскольку объединение $i$-ых лучей в двух разрезах, разбивает плоскость на две части, одна из которых полностью содержится в $\alpha_i$, а, значит, не пересекается с объединением остальных $n-1$ лучей.
\medskip

Верен аналог теоремы~\ref{pr:bij-2} 
для плоскости без $n$ точек.

\begin{theorem}\label{p:bij-n}
Следующее отображение является взаимно однозначным соответствием между множеством гомотопических классов замкнутых ломаных на плоскости без $n$ точек и множеством $E_{n,c}$ экономных циклических слов из букв $p_1,\ldots,p_n$, $p_1^{-1},\ldots,p_n^{-1}$.
Отображение переводит замкнутую ломаную $l$ в экономную форму $e(l)$ ее циклического слова Пуанкаре относительно некоторого разреза общего положения, чьи лучи лежат в $\alpha_1, \ldots, \alpha_n$, соответственно.
\end{theorem}


По утверждению \ref{p:invar-eco-n}.b, отображение из формулировки теоремы \emph{корректно определено} на множестве гомотопических классов.
Аналогично примеру \ref{p:cycpoi}, это отображение \emph{сюръективно}.
Остается доказать его \emph{инъективность}.

\begin{proof}[Набросок доказательства инъективности при $n\geqslant2$]
(Ср. с наброском доказательства утверждения~\ref{veplho-poly}.b.)
Выберем точку $O$, не лежащую на одной прямой ни с какими двумя из точек $P_1, P_2, \ldots, P_n$.
Выберем $n$ лучей с началом в точке $O$, так чтобы их объединение разбивало плоскость на $n$ внутренностей (выпуклых) углов $\beta_1, \ldots, \beta_n$, содержащих точки $P_1,  \ldots, P_n$, соответственно.
Внутри каждого угла $\beta_i$ выберем две точки $S_i$ и $T_i$, так чтобы лучи $P_iS_i$ и $P_iT_i$ лежали внутри $\beta_i$ и для треугольника $\Delta_i := OS_iT_i$ было верно $w(\Delta_i, P_i) = +1$.

Достаточно доказать, что \textit{любая замкнутая ломаная гомотопна замкнутой ломаной,} 

$(*)$ \textit{каждое невырожденное звено которой является стороной одного из треугольников $\Delta_1, \ldots, \Delta_n$}. 

Докажем это.
Назовем пару $\{ XY,\ B \}$, состоящую из отрезка $XY$ ломаной $l$ и вершины $B \notin \{X,Y\}$ одного из треугольников $\Delta_1, \ldots, \Delta_n$, \emph{неудобной} (для $l$), если либо 

$\bullet$ $O\ne B$ и $B \in \angle XP_iY$, либо  

$\bullet$ $O=B$ и отрезок $XY$ пересекает объединение $OP_1 \cup OP_2 \cup \ldots \cup OP_n$.

Назовем ломаную \emph{удобной}, если 

(1) она содержит точку $O$, и каждый участок ломаной между <<соседними>> вхождениями точки $O$ целиком содержится в одном из углов $\beta_1, \ldots, \beta_n$, и 

(2) ни для каких звена $XY$ ломаной $l$ и вершины $B$ одного из треугольников $\Delta_1, \ldots, \Delta_n$ пара $\{XY, B \}$ не является неудобной.
   
В следующем абзаце мы доказываем, что \emph{любая замкнутая ломаная $l$ гомотопна удобной.} 

Можно считать, что $l$ не является удобной. 
Сделав несколько (если необходимо) обратных операций к элементарному сокращению, мы можем считать, что $l$ удовлетворяет свойству~(1).
Тогда существует неудобная пара $\{ XY,\ B \}$ для $l$.
Добавим $B$ к $l$ между $X$ и $Y$, используя обратную операцию к элементарному сокращению.
Полученная ломаная удовлетворяет свойству~(1).
Количество неудобных пар для полученной замкнутой ломаной меньше, чем для $l$.
С помощью конечного числа таких операций получим удобную замкнутую ломаную, гомотопную $l$.

Любую вершину удобной замкнутой ломаной, отличную от вершин треугольников $\Delta_1, \ldots, \Delta_n$, можно удалить элементарным сокращением.
Полученная замкнутая ломаная также удобна.
Следовательно, с помощью конечного числа таких элементарных сокращений мы получим удобную замкнутую ломаную, содержащую только вершины треугольников $\Delta_1, \ldots, \Delta_n$.
Тогда для полученной ломаной выполнено $(*)$.
\end{proof}

\aronly{
\noindent\emph{Указание к другому доказательству инъективности (А. Абзалилов).}
Для каждого экономного циклического слова $\alpha$ из букв $p_1,\ldots,p_n$, $p_1^{-1},\ldots,p_n^{-1}$ построим замкнутую ломаную $k(\alpha)$ и докажем, что $l$ гомотопна $k(e(l))$.
Выберем точку $X$, не лежащую на одной прямой ни с какими двумя из точек $P_1, P_2, \ldots, P_n$.
Отметим на плоскости точки $S_1$ и $T_1$ так, чтобы $S_1, T_1$ и $X$ не лежали на одной прямой, внутри треугольника $S_1T_1X$ из выколотых лежала только точка $P_1$, и $w(XS_1T_1, P_1)=+1$.
Аналогично определим точки $S_i$ и $T_i$ для каждого $i=1,\ldots,n$.
Есть только одна замкнутая ломаная $k(\alpha)$, соответствующая $\alpha$, звенья которой --- стороны треугольников $S_iT_iX$ и никакие два соседних ребра не направлены в противоположные стороны.

Построим гомотопию, переводящую $l$ в $k(e(l))$.
Обозначим соседей вершины $B$ ломаной $l$ через $B_-$ и $B_+$.

Пусть некоторая вершина $B$ текущей ломаной лежит на одном из лучей $P_iP_i'$.
Отметим на отрезке $BB_+$ точку $B'$ так, чтобы треугольник $BB'B_-$ не содержал выколотых точек.
Заменим $B_-BB'$ на $B_-B'$ элементарным сокращением.
После таких операций \emph{ни одна вершина текущей ломаной не будет лежать ни на одном из лучей $P_iP_i'$}.

\emph{Удлинением} назовем операцию, обратную элементарному сокращению.

Пусть некоторое звено $B_-B$ текущей ломаной пересекает несколько лучей $P_iP_i'$.
Добавим удлинением точки на отрезке $B_-B$, чтобы никакие из получившихся звеньев не пересекали более одного луча $P_iP_i'$.
После таких операций \emph{каждое звено текущей ломаной пересекает не более одного 
луча $P_iP_i'$}.

\medskip
}





\begin{figure}[h]\centering
\compile{
\includegraphics[scale=1.]{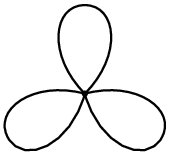} 
}
\caption{Граф <<букет трех окружностей>>}\label{f:wedge}
\end{figure}

Верен аналог загадки \ref{p:defret} 
с заменой плоскости без двух точек на плоскость без $n$ точек, а графа <<восьмерка>> "--- на граф <<букет $n$ окружностей>> (рисунок~\ref{f:wedge}).

\section{Умножение замкнутых ломаных}\label{s:mulpol}

\textbf{Отмеченной замкнутой} ориентированной \textbf{ломаной} $A_1\ldots A_m$ на плоскости называется последовательность (т.~е. упорядоченный набор) $(A_1,\ldots,A_m)$ точек плоскости (не обязательно различных).
Далее слова <<ориентированная>> и <<на плоскости>> пропускаются. 
Точка $A_1$ называется \emph{отмеченной точкой}.


\textbf{Обратной} к отмеченной замкнутой ломаной $l = A_1 A_2 \ldots A_m$ называется отмеченная замкнутая ломаная $l^{-1} := A_1 A_m \ldots A_2$.

В этом параграфе $l_1, l_2$ "--- отмеченные замкнутые ломаные с общей отмеченной точкой $X$.

\begin{figure}[h]\centering
\compile{
\includegraphics[scale=0.3]{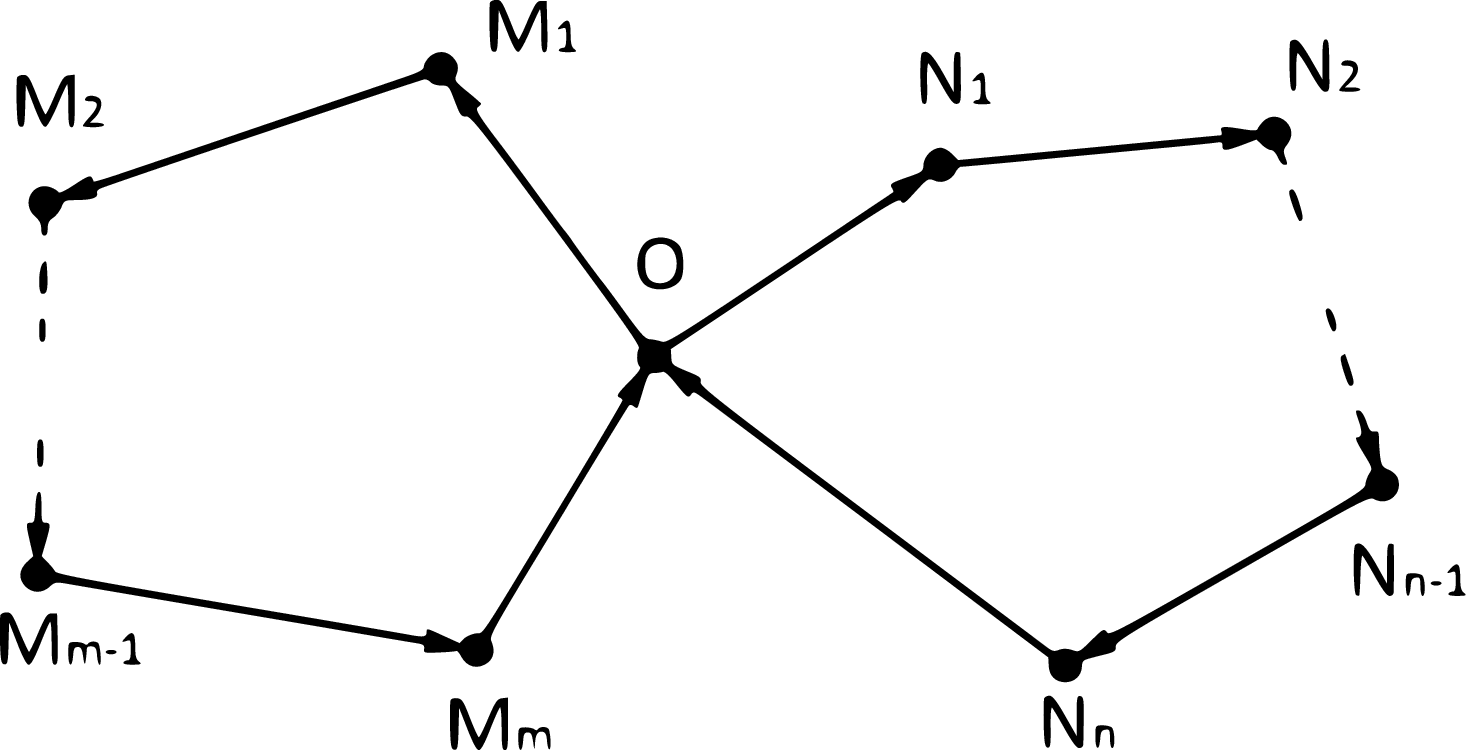}
}
\caption{Произведение отмеченных замкнутых ломаных $l_1 = OM_1\ldots M_m$ и $l_2 = ON_1\ldots N_n$}
\label{f:konk}
\end{figure}

\textbf{Произведением} (конкатенацией, соединением) отмеченных замкнутых ломаных  $l_1 = XM_1\ldots M_m$ и $l_2 = XN_1\ldots N_n$ с общей отмеченной точкой называется отмеченная замкнутая ломаная
$l_1l_2 := XM_1 \ldots M_m X N_1 \ldots N_n$ (см. рисунок~\ref{f:konk}).

\begin{remark*}
На плоскости без двух точек $P,Q$ возьмем отмеченные замкнутые ломаные $a$ и $b$ (треугольники на рисунке \ref{f:eight}) с общей отмеченной точки, ограниченные которыми части плоскости пересекают $\{P,Q\}$ по $P$ и по $Q$, соответственно. 
В этом замечании мы сокращаем <<замкнутая ломаная, полученная из отмеченной замкнутой ломаной $x$ забыванием отмеченной точки>> до <<замкнутая ломаная $x$>>. 

Примером \ref{veplho-borr}.a является замкнутая ломаная $[a, b] := aba^{-1}b^{-1}$. 
Примером \ref{veplho-borr}.b будет $[a, b]c[a, b]^{-1}c^{-1}$ (сделайте это строгим самостоятельно).

Вот версия парадокса Пуанкаре.
\emph{Замкнутые ломаные $aba^{-1}$ и $b$ гомотопны, но замкнутые ломаные $aba^{-1}b^{-1}$ и $bb^{-1}$ не гомотопны (ломаная $bb^{-1}$ гомотопна одноточечной замкнутой ломаной).}
\end{remark*}

\begin{proposition}[аддитивность]
\label{p:add} 
Число оборотов произведения $l_1l_2$ вокруг точки $O$ вне $l_1,l_2$ равно сумме чисел оборотов ломаных $l_1$ и $l_2$ вокруг $O$.
\end{proposition}

\begin{proposition}
\label{p:ngroup} 
Умножение отмеченных замкнутых ломаных с общей отмеченной точкой $X$

(a) ассоциативно, но не имеет единичного элемента; 

(b) (загадка) не порождает корректно определенного умножения на классах гомотопности.
\end{proposition}

Отмеченные замкнутые ломаные с общей отмеченной точкой называются \emph{отмеченно гомотопными}, если они <<гомотопны с сохранением отмеченной точки>> (сделайте это строгим самостоятельно).

\begin{proposition}
\label{p:based-nonbased} 
Если отмеченные ломаные гомотопны в плоскости без точки $O$, то они также отмеченно гомотопны в плоскости без точки $O$.
\end{proposition}

\begin{example}
Существуют отмеченные ломаные, гомотопные в плоскости без двух точек, не являющиеся отмеченно гомотопными в плоскости без двух точек. 
\end{example}

\begin{proposition}
\label{p:group} 
(a) Умножение классов отмеченной гомотопности отмеченных замкнутых ломаных корректно определено, ассоциативно и имеет единичный элемент.

(b) Для этого умножения каждый класс имеет обратный.
\end{proposition}

\begin{proof}[Набросок доказательства корректности из утверждения~\ref{p:group}.a] (получен редактированием текста Н. Буровик)
Обозначим через $l_1, l'_1, l_2, l'_2$ отмеченные замкнутые ломаные с общей отмеченной точкой. 
Пусть $l_1, x_1, x_2, \ldots, x_n, l'_1$ и $l_2, y_1, y_2, \ldots, y_m, l'_2$ "--- отмеченные гомотопии (между $l_1$ и $l'_1$, и между $l_2$ и $l'_2$). 
Тогда 
$l_1l_2, x_1l_2, x_2l_2, \ldots, x_nl_2, l'_1l_2, l'_1y_1, l'_1 y_2, \ldots, l'_1 y_m, l'_1l'_2$ "--- отмеченная гомотопия (между $l_1 l_2$ и $l'_1 l'_2$).
\end{proof}

Таким образом, естественная групповая структура существует не на интересующем нас множестве (множество замкнутых ломаных с точностью до гомотопности), а на множестве, которое в данный момент кажется менее естественным (множество отмеченных замкнутых ломаных с точностью до отмеченной гомотопности).
Действительно, чтобы определить произведение, нужны \emph{отмеченные} замкнутые ломаные; 
чтобы получить корректно определенное умножение на гомотопических классах, нужна \emph{отмеченная} гомотопность (утверждения \ref{p:ngroup}.b и \ref{p:group}).

Далее, отмеченные замкнутые ломаные и их отмеченная гомотопия рассматриваются \emph {на плоскости без двух точек $P$ и $Q$}.

\emph{Слово Пуанкаре} отмеченной замкнутой ломаной $x$ определяется аналогично циклическому слову Пуанкаре, только движение вдоль $x$ начинается и заканчивается в отмеченной точке.

\begin{example}[ср. с примером~\ref{p:cycpoi}]
\label{p:surj-2-based}
Любое слово из букв $a,b,a^{-1},b^{-1}$ является словом Пуанкаре некоторой отмеченной замкнутой ломаной.
\end{example}

Аналогично~\S\ref{s:minus2}, можно определить экономную форму слова из букв $a,b,a^{-1},b^{-1}$.
Экономная форма слова Пуанкаре отмеченной замкнутой ломаной является инвариантом отмеченной гомотопности (аналогично утверждениям~\ref{p:invar-2eco}.b,~\ref{p:invar-eco}.b). 
Далее мы приводим другое определение инварианта отмеченной гомотопности в терминах отношения эквивалентности на множестве всех слов из букв $a,b,a^{-1},b^{-1}$ (аналогичное замечанию \ref{r:equiv}). 

Рассмотрим множество всех (конечных) слов (включая пустое слово) из букв $a, b, a^{-1}, b^{-1}$.
Для такого слова \emph{элементарное сокращение} — это замена любого из подслов $aa^{-1}, a^{-1}a, bb^{-1}, b^{-1}b$ пустым подсловом.
Два таких слова называются \textit{эквивалентными}, если их можно соединить последовательностью слов, в которой одно из любых двух последовательных слов может быть получено из другого с помощью элементарного сокращения.
Обозначим через $F_2$ множество всех классов эквивалентности (стандартное обозначение: $\left<a,b\right>$).

Класс эквивалентности слова Пуанкаре отмеченной замкнутой ломаной может измениться  при изменении лучей $PP'$ и $QQ'$. 
(Это аналогично примеру~\ref{e:3raysn}.)

\begin{proposition}[ср. с утверждением~\ref{p:invar-2eq}]\label{p:invar-eq-based}
Класс эквивалентности слова Пуанкаре не меняется при замене отмеченной замкнутой ломаной на отмеченно гомотопную ей.
\end{proposition}

\begin{proposition}[ср. с утверждением~\ref{p:comp-eco}]
\label{p:comp-eq-based}
Если слова Пуанкаре двух отмеченных замкнутых ломаных $l_1, l_2$ эквивалентны, то $l_1$ и $l_2$ отмеченно гомотопны.
\end{proposition}

Умножение слов порождает операцию на $F_2$ и превращает $F_2$ в группу (эта группа называется \emph{свободной группой с двумя образующими}).
Поэтому следующая теорема~\ref{pr:bij-2-based} следует из примера~\ref{p:surj-2-based} и утверждений~\ref{p:invar-eq-based}, \ref{p:comp-eq-based}.

\begin{theorem}[ср. с теоремой~\ref{pr:bij-2}]
\label{pr:bij-2-based} 
Следующее отображение является изоморфизмом между группой классов отмеченной гомотопности отмеченных замкнутых ломаных на плоскости без двух точек и группой $F_2$.
Отображение переводит отмеченную замкнутую ломаную в ее слово Пуанкаре.
\end{theorem}

\begin{remark*}
Заметим, что $F_{2,c}$ (см. определение в замечании после теоремы~\ref{pr:bij-2}) "--- множество классов \emph{сопряженности} группы $F_2$.
\end{remark*}

\begin{theorem}
\label{hompi-alg} 
(a)~Для любого графа $G$ с отмеченной точкой существует алгоритм, распознающий отмеченную гомотопность отмеченных ориентированных циклов в $G$. 

(b) Для любого графа $G$ существует алгоритм, распознающий гомотопность ориентированных циклов в~$G$. 

(c) Для любых двух графов существует алгоритм, распознающий гомотопность \emph{симплициальных отображений} между ними.
(См. определение в \cite[\S\,9 <<Гомотопическая классификация отображений>>]{Sk}, \cite{Sk20e}.)
\end{theorem}

Аналоги этой теоремы для $2$-гиперграфов неверны! 
Это следует из алгоритмической неразрешимости проблемы тривиальности в некоторой группе, определяемой конечным числом образующих и соотношений. 
См.~\cite[теорему 14.3.1]{Sk20} и ее обобщение в~\cite{Sk20e}.

{\it В этом списке звездочками отмечены книги, обзоры и популярные статьи.}

%% file: obstr-hoclpoli_eng.tex
\section{Introduction}

In this text we expose (as a sequence of problems) basic cases of some fundamental ideas and methods of mathematics. 
Namely, of homotopy, degree, fundamental group, covering, Whitehead invariant, etc.
This is done by considering the elementary example: closed polygonal lines in a subset of the plane. 
Although these ideas and methods are parts of topology, they are used in many other areas including computer science~\cite{DDM+}. 
There are both simple (\ref{p:veplho-poly0b}), and complicated (\ref{3space}) problems.

This text is expository and is accessible to mathematicians not specialized in the area (and to students).

An \textbf{oriented cyclic sequence} is an ordered set up to a cyclic shift. 
A \textbf{closed} oriented \textbf{polygonal line} in the plane is an oriented cyclic sequence of points in the plane (the points need not be distinct).
Below, the words `oriented' and `in the plane' are omitted.\footnote{Thus, a closed polygonal line (defined here) is not a subset of the plane. 
Still, we sometimes work with a closed polygonal line $A_1\ldots A_m$ as with the union of segments $A_iA_{i+1}$, e.g. we write `a closed polygonal line, not passing through a point'.}

\begin{figure}[h]\centering
\includegraphics[scale=0.5]{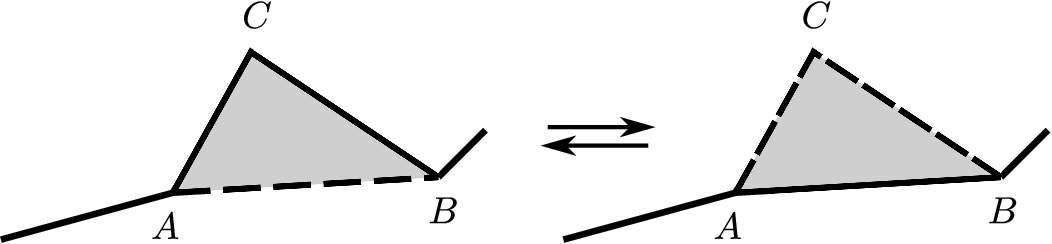}
\caption{An elementary cancellation}
\label{f:elem}
\end{figure}

Let $N$ be a subset of the plane. 
Informally speaking, two closed polygonal lines in~$N$ are said to be \emph{homotopic}, if one can be transformed to the other by a `continuous deformation' in $N$.
Let us give a rigorous definition.
An \emph{elementary cancellation} in~$N$ of a closed polygonal line in~$N$ is the removal of a vertex~$B$ such that for vertices $A,C$ (that need not be distinct) adjacent to~$B$ the convex hull of vertices $A,B,C$ \emph{is contained in~$N$} (see Figure~\ref{f:elem}).
A (piecewise linear) \textbf{homotopy in $N$} is a finite sequence of closed polygonal lines in which for any two consecutive polygonal lines one is obtained from the other by an elementary cancellation in~$N$.
Two closed polygonal lines are said to be (piecewise linearly) \textbf{homotopic} in~$N$, if there exists a homotopy in $N$, whose first and last polygonal lines coincide with the given ones.

\emph{Below, when we write about a closed polygonal line in some subset of the plane, we consider the property of being homotopic (and elementary cancellations) in this very subset.}

\begin{pr}\label{p:veplho-poly0b}
Any two one-point closed polygonal lines in the plane minus a point are homotopic.
\end{pr}

\begin{proof}
Given points~$A$ and~$B$, take a point~$C$ not lying on any of the lines joining $A$ and $B$ to the deleted point.
Then $A,AC,C,CB,B$ is a required homotopy.
\end{proof}

The homotopy classification of closed polygonal lines \emph{in the plane} is trivial.

\begin{pr}\label{p:veplho-poly0a} Any two closed polygonal lines in the plane are homotopic.
\end{pr}

\begin{figure}[h]\centering
\includegraphics[scale=0.7]{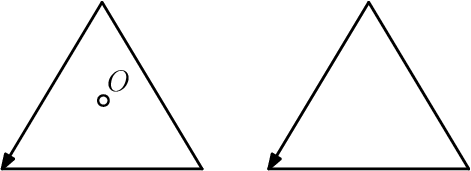}
\caption{Two closed polygonal lines in the plane minus a point $O$}
\label{f:guess}
\end{figure}

\begin{pr}\label{p:veplho-poly1} 
There are non-homotopic closed polygonal lines in the plane minus a point.
\end{pr}

\begin{figure}[H]
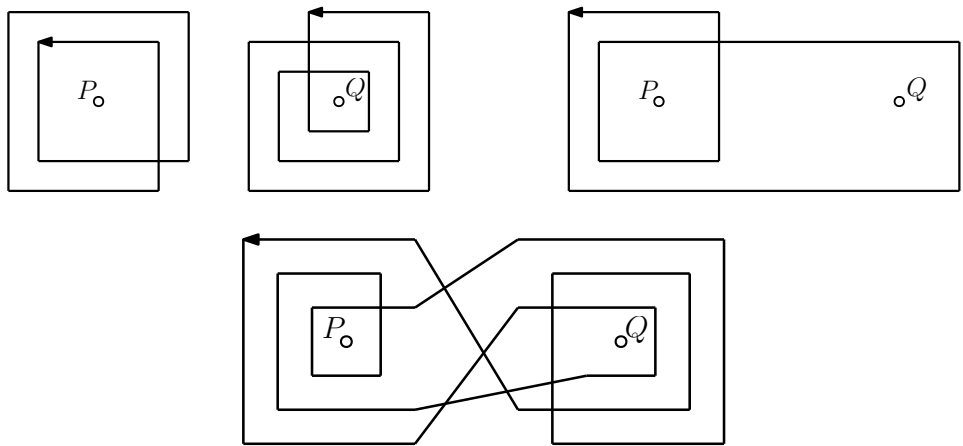

\centering
\includegraphics[scale=0.7]{2points/a2b3.eps} \qquad \qquad
\includegraphics[scale=0.7]{2points/aab.eps}\\
\vspace{0.5cm}
\includegraphics[scale=0.8]{2points/aab-abb.eps}
\caption{The four closed polygonal lines in the plane minus two points $P, Q$}
\label{f:2points-a2-b3}
\end{figure}

\begin{figure}[H]
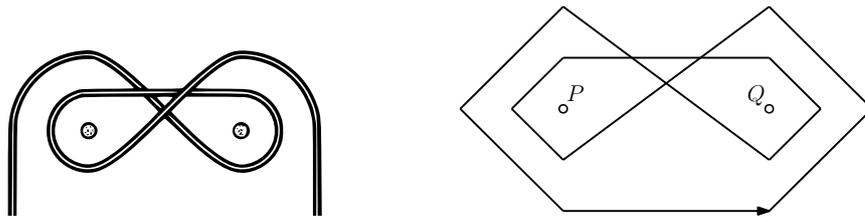
\centering
\compile{
\includegraphics[scale=0.25]{two-nails.eps}
\qquad
\qquad
\includegraphics[scale=0.6]{2points/poincare.eps}
}
\caption{To Problem~\ref{veplho-borr}.a: a rope hanging on two nails (left); a closed polygonal line in the plane minus two points $P, Q$ (right)}
\label{f:poincare}
\end{figure}

\begin{figure}[H]\centering
\compile{
\includegraphics[scale=0.9]{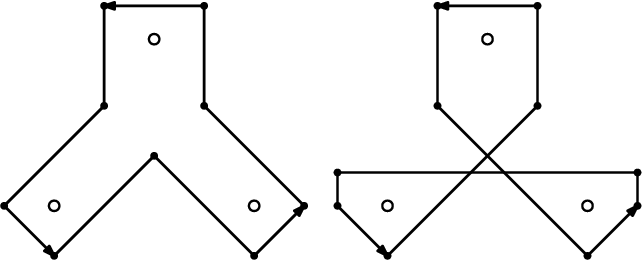}
}
\caption{Two closed polygonal line in the plane minus three points}
\label{f:3points}
\end{figure}

See e.g.~Figure~\ref{f:guess}. 
Try to invent invariants that distinguish the closed polygonal lines of Problems~\ref{p:veplho-poly1}--\ref{p:3points} and~\ref{3space} (in particular, in Figures~\ref{f:guess}--\ref{f:poincare} and \ref{f:3points}).
The definition of the invariants, and their properties form the theory presented in~\S\S\ref{s:wind}--\ref{s:minus2}.
This theory is also required for Theorem~\ref{t:hocl-algor}.

\begin{pr}\label{veplho-poin} 
The four closed polygonal lines in Figure~\ref{f:2points-a2-b3} are pairwise not homotopic in the plane minus two points.
\end{pr}

\begin{pr}[the Poincar\'e paradox]\label{veplho-borr} (a) How can one hang a closed rope (with a heavy medal) around two nails (driven into a flat wall) so that the rope does not fall, but removing either nail causes the rope to fall?

More rigorously, give an example of two points $P,Q$ and a closed polygonal line in~$\R^2-P-Q$, that is homotopic to a one-point closed polygonal line in~$\R^2-P$, and in~$\R^2-Q$, but not in~$\R^2-P-Q$.

(b) An analogous problem for three nails.
\end{pr}

These examples (and the example in \ref{3space}.a) can first be given without proof of their properties.

\begin{figure}[H]\centering
\compile{
\includegraphics[scale=0.3]{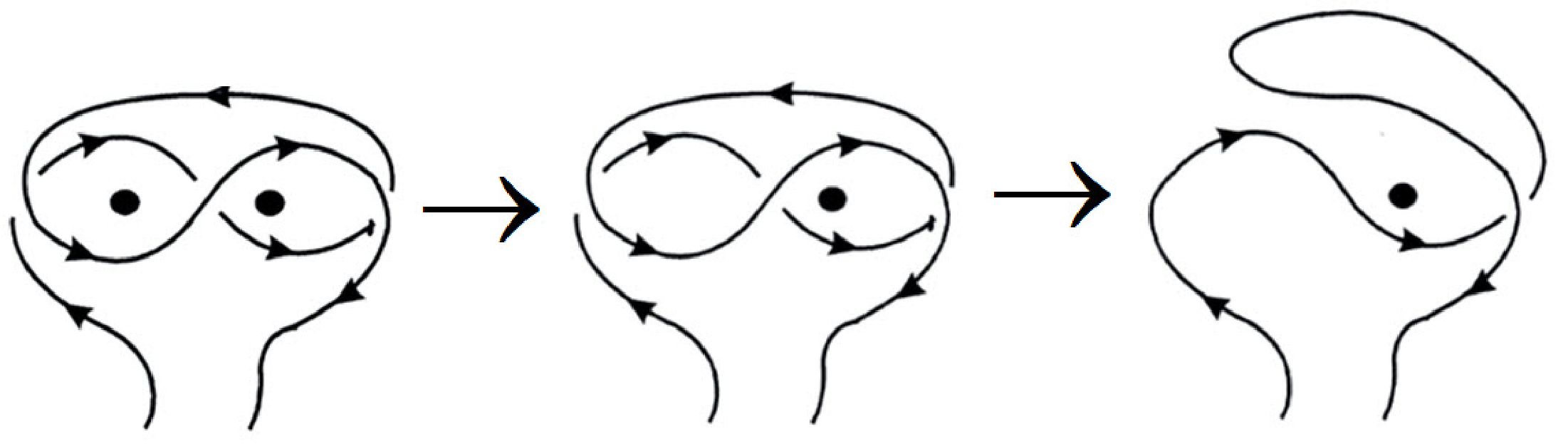}
}
\caption{Removing either nail causes the rope to fall}
\label{f:drop}
\end{figure}

\begin{figure}[H]\centering
\compile{
\includegraphics[scale=0.4]{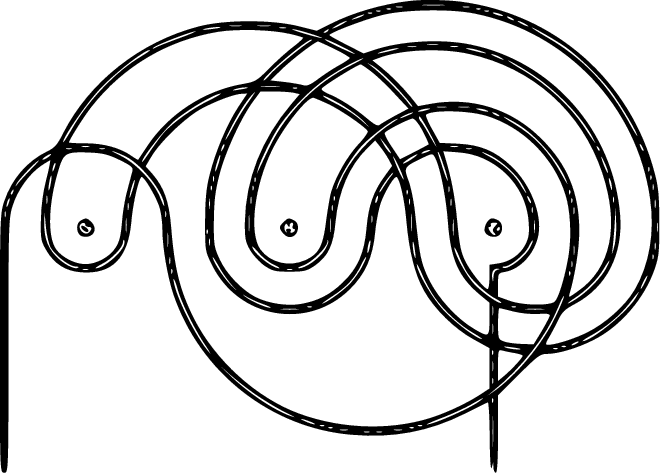}
}
\caption{To Problem~\ref{veplho-borr}.b: a rope hanging on three nails}
\label{f:three}
\end{figure}

\smallskip
Examples to Problem~\ref{veplho-borr} are given in Figures~\ref{f:poincare}--\ref{f:three}, or in the remark at the beginning of \S\ref{s:mulpol}.
See also~\cite{Zi10},~\cite{DDM+}.

\begin{pr}\label{p:3points}*
The closed polygonal lines in Figure~\ref{f:3points} are not homotopic in the plane minus three points.
\end{pr}

By studying this text, the reader will learn useful and important ideas and methods. 
They are useful and important because they give the following bright result, whose statement is accessible to a non-specialist in topology (and many other results, see textbooks on topology, e.g.~\cite{Sk20}).

\begin{theorem}\label{t:hocl-algor} For any $n$, there is an algorithm recognizing homotopy of closed polygonal lines in the plane minus $n$ points.
\end{theorem}

This is an endpoint of the main results on homotopy classification in the plane minus $n$ points (starting from $n=1,2$, see Propositions~\ref{pr:bij-1},~\ref{pr:bij-2}, and 
Problem~\ref{p:bij-n}.a).

You will see how homotopy classification is related to \emph{combinatorics of words}.
A group structure will naturally appear, although on a slightly different set of \emph{based} closed polygonal lines up to \emph{based} homotopy.
This illuminates the relation of the exposed classification to \emph{combinatorial group theory}. 
Due to this relation, topological methods can be used in group theory. 
This is the starting point of \emph{geometric group theory}~\cite{GG}.   

\begin{figure}[h]
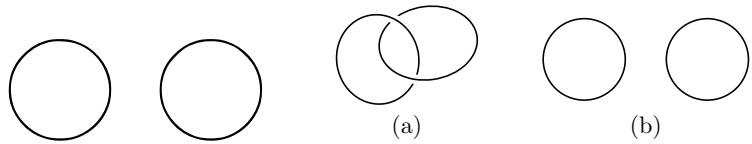
\centering
\compile{
\includegraphics[scale=0.25]{2_rings.eps} \qquad \includegraphics[scale=0.9]{3-1.eps}
}
\caption{Unlinked and linked rings}\label{f:hopf}
\end{figure}

\begin{pr}[Riddle]\label{3space}* In 3-space there are two 

(a) unlinked; \quad (b) linked  

rings as in Figure~\ref{f:hopf}. 
Can one wind and close a rope so that the closed rope cannot be pulled far away from the two rings, but can be pulled far away from one ring after cutting the other ring? 

The formalization is analogous to Problem \ref{veplho-borr}. 
\end{pr}

\medskip
\section{On the style of this text}
\emph{In this text we expose a theory as a sequence of problems,} see e.g. \cite{HC19}, \cite[Introduction, Learning by doing problems]{Sk21m} and the references therein. 
Most problems are useful theoretical facts. 
So this text could in principle be read even without solving problems.
If a mathematical statement is formulated as a problem, then the objective is to prove this statement.
Open-ended questions are called {\bf riddles}; here one must come up with a clear wording, and a proof. 
\emph{If a problem is named `theorem' (`lemma', `corollary', etc.), then this statement is considered to be more important.}
Usually we \emph{formulate} beautiful or important statements \emph{before} giving a sequence of results (lemmas, propositions, etc.) which constitute its \emph{proof}.
We give hints on that after the statements but we do not want to deprive you of the pleasure of finding the right moment when you finally are ready to prove the statement.
In general, if you are stuck on a certain problem, try looking at the next ones; they may turn out to be helpful.
Problems marked by star and remarks are not used in the sequel; although problems with star are not necessarily complicated, they can be omitted during the first round of problem solving.
Important definitions are highlighted in \textbf{bold} for easy navigation.

\section{Winding number: definition and discussion}\label{s:wind}

In the rest of this text, all points and closed polygonal lines are considered in the plane.

Let $O, A, B$ be points such that $A\ne O$ and $B\ne O$ (but possibly $A=B$). 
The \emph{oriented (a.k.a. directed) angle} $\angle AOB$ is the number $t\in(-\pi,\pi]$ such that the vector $\overrightarrow{OB}$ is codirected to the vector obtained from $\overrightarrow{OA}$ by the rotation through $t$.   
(If you can treat vectors in the plane as complex numbers, then you can rewrite this condition as $\overrightarrow{OB}\upuparrows e^{it}\overrightarrow{OA}$.)
Below, oriented angles are considered, and the word `oriented' is omitted.

When solving problems, the following statement (close to axioms) can be used without proof:
{\it For any points $A,B,C$ in the plane and a point $O$ not lying on the union of segments $AB,BC,CA$
 
$\bullet$ $\angle OAB + \angle OBC + \angle OCA = \pm 2\pi$, if $O$ is in the convex hull of points $A,B,C$,

$\bullet$ $\angle OAB + \angle OBC + \angle OCA = 0$, otherwise.}

Let $l = A_1\ldots A_m$ be a closed polygonal line not passing through $O$. 
The {\bf winding number} $w(l,O)$ of $l$ around $O$ is said to be the number of revolutions during the rotation of the vector whose origin is $O$, and whose endpoint goes along the polygonal line in a positive direction.
Rigorously, 
$$2\pi \cdot w(l,O) := 
\angle A_1OA_2+\angle A_2OA_3+\ldots+\angle A_{m-1}OA_m+\angle A_mOA_1.$$

\begin{figure}[ht]
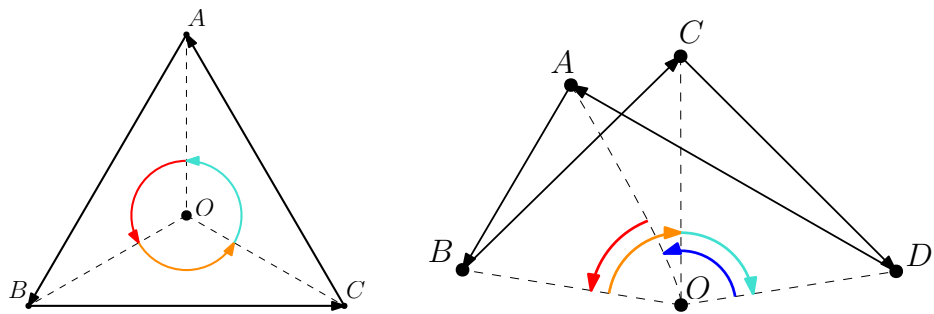
\centering
\compile{
\includegraphics[scale=0.65]{abco.eps}\qquad
\includegraphics[scale=0.85]{all_angles.eps}
}
\caption{$w(ABC, O) = +1 $ and $w(ABCD,O) = 0$}
\label{f:abco}
\end{figure}
 
E.g. in Figure~\ref{f:abco}
$w(ABC, O) = \dfrac{1}{2\pi} \left(\angle AOB + \angle BOC + \angle COA \right) = +1\quad\text{and}$
$$
2\pi \cdot w(ABCD,O)=\angle AOB + \angle BOC + \angle COD + \angle DOA = \angle BOD + \angle DOB = 0.
$$

\begin{figure}[ht]\centering
\compile{
\includegraphics[scale=0.9]{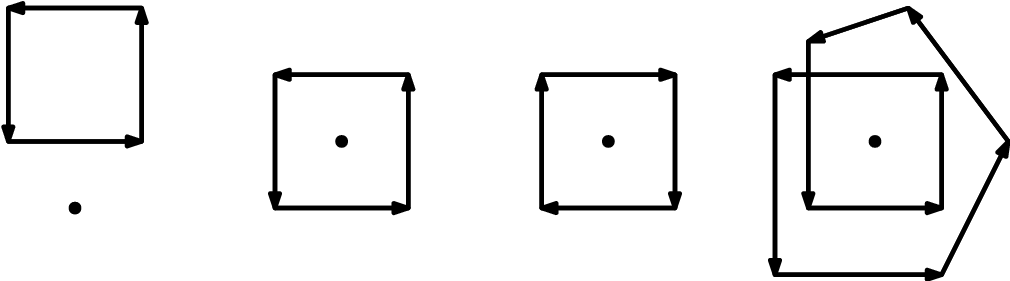}
}
\caption{The winding numbers equal $0,~+1,~-1,~+2$}
\label{f:wn_examples}
\end{figure}
 
\begin{pr}\label{p:winsimdeg} (a) Let $ABC$ be a regular triangle and $O$ its center. 
Find $w(ABCABC,O)$. 

(b) Give an example of a closed polygonal line $l$ in the plane such that $w(l,O)=0$ for any point $O\in\R^2-l$.
\end{pr}

The result of Problem~\ref{p:winsimdeg}.a shows that winding numbers for distinct closed polygonal lines with the same union of their segments can be distinct.

\begin{proposition}\label{p:winsim} The winding number of the outline of any convex polygon around any point in its exterior (interior) is $0$ ($\pm 1$).
See Figure~\ref{f:wn_examples}.
\end{proposition}

\begin{pr}\label{p:surj-wn} For any integer $n$ and any point $O$, there is a closed polygonal line whose winding number around $O$ is $n$. 
\end{pr}

\begin{proposition}\label{p:noncl} The winding number $w(A_1\ldots A_m,O)$ is an integer.
\end{proposition}

\emph{Hint.} By the bullet point properties at the beginning of this section, 
$$\angle A_{m-1}OA_m+\angle A_mOA_1 \equiv \angle A_{m-1}OA_1 \mod2\pi.$$  

See hints and solutions in \cite[\S1]{ABM+r}, \cite[\S1, \S8]{ABM+e}. 
For more on winding number and related notions see~\cite{ABM+r}, \cite{ABM+e} and the references therein.

\section{Homotopy classification in the plane minus a point}

In this section, closed polygonal lines and their homotopy are considered \emph{in the plane minus a point $O$}.

\begin{pr}\label{veplho-poly}
(a)~Closed polygonal lines with distinct winding numbers around~$O$ are not homotopic.

(b)~A closed polygonal line $l$ is homotopic to a one-point closed polygonal line if and only if $w(l,O)=0$.

(c)~Take the outline $\partial\Delta$ of a triangle $\Delta$ such that a point~$O$ is inside the triangle.
Any closed polygonal line is homotopic to $|w(l,O)|$-times 
traversal of the outline $\partial\Delta$, the traversal made 

$\bullet$ counterclockwise if $w(l,O)\ge0$,

$\bullet$ clockwise if $w(l,O)\le0$.

(Then closed polygonal lines having the same winding numbers around~$O$ are homotopic.)
\end{pr}

\begin{proposition}\label{pr:bij-1} 
The following map between the set of homotopy classes of closed polygonal lines in the plane minus a point, and the set of integers is a 1--1 correspondence. 
The image of a closed polygonal line $l$ (more precisely, of its homotopy class) is the winding number of $l$ around the point.

(Proposition~\ref{pr:bij-1} is proved in Problems~\ref{p:surj-wn}~and~\ref{veplho-poly}.)
\end{proposition}

A ray in the plane is said to be \emph{in general position w.r.t. a closed polygonal line $l$} if the ray does not pass through the vertices of $l$.

\begin{pr}\label{p:ray} * Take a closed polygonal line $l$, and a ray $OP$ in general position  w.r.t.~$l$. 

(a) Does the number of intersection points of $l$ and $OP$ have the same parity as $w(l,O)$?
    
(b) The number of segments of $l$ that intersect $OP$ has the same parity as $w(l,O)$.
\end{pr}

\emph{Hint.} For a proof of (b) using homotopy, 
Problems~\ref{veplho-poly}.c and \ref{p:apprgen} are useful.
For a proof of (b) not using homotopy, see \cite[\S2]{ABM+r}, \cite[\S2]{ABM+e}.

\begin{pr}\label{p:apprgen}*
Take two homotopic closed polygonal lines $l, m$.
Take a ray $OP$ in general position w.r.t. $l$ and $m$.
Then there is a homotopy between $l$ and $m$ such that the ray $OP$ is in general position w.r.t. every polygonal line of the homotopy.
\end{pr}

It is useful to understand that the `2-dimensional' Proposition~\ref{pr:bij-1} is in essence `1-dimensional'. 
Analogous understanding is used for homotopy classification in the plane without more than one point. 

\begin{figure}[H]\centering
\compile{
\includegraphics[scale=0.25]{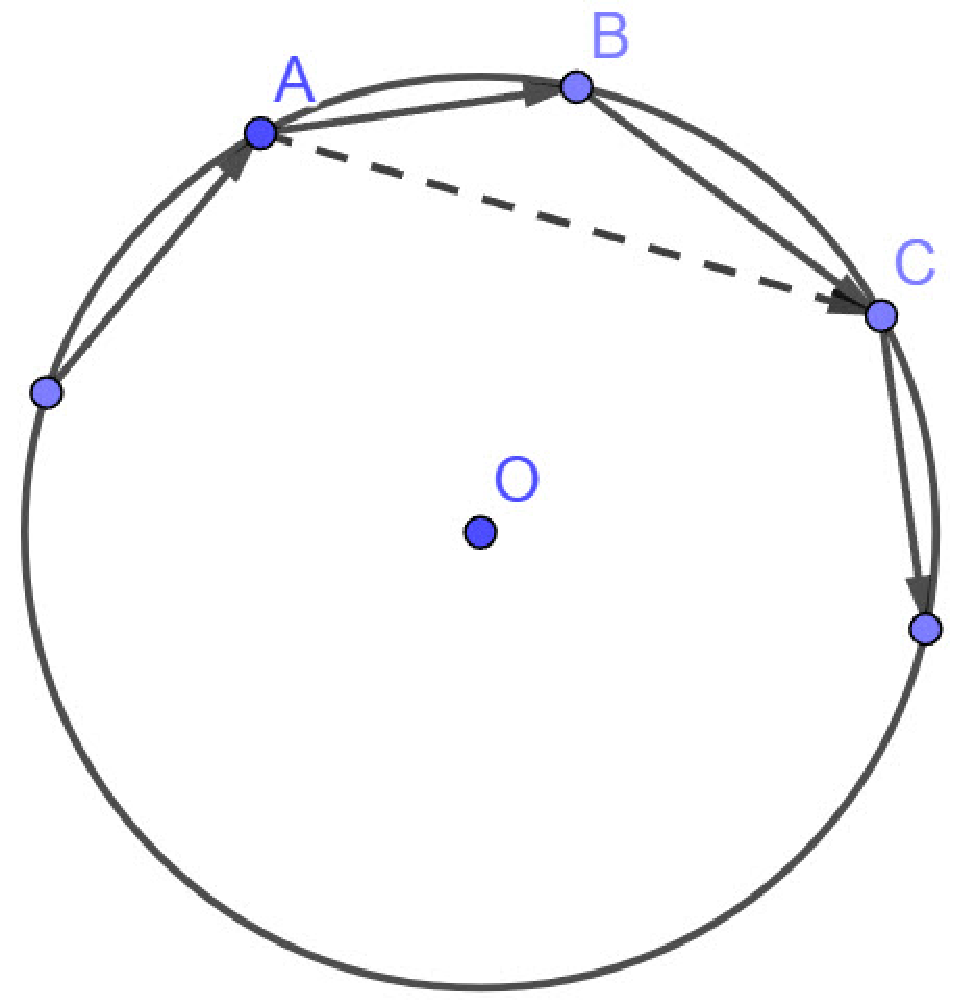}
}
\caption{An elementary cancellation of a traversal}
\label{f:rou}
\end{figure}

\begin{pr}\label{veplho-polyc}* 
A \emph{traversal} is an oriented cyclic sequence of points on a (fixed) circle such that no two consecutive points are diametrically opposite.
An \emph{elementary cancellation} of a traversal is the removal of a point~$B$ such that for points $A,C$ (which need not be distinct) adjacent to~$B$ the convex hull of points $A,B,C$ does not contain the center $O$ of the circle (see Figure~\ref{f:rou}).  
Two traversals are said to be \emph{homotopic}, if there is a sequence of traversals, in which one of any two consecutive traversals can be obtained from the other by an elementary cancellation.

(a)~A traversal~$l$ is homotopic to a one-point traversal if and only if $w(l,O)=0$. 

(b)~(Riddle) Describe traversals up to homotopy.
\end{pr}

\begin{proof}[Hints to Problem~\ref{veplho-poly}]
    (a) It is sufficient to prove that \textit{if a closed polygonal line~$l_1$ is obtained from a closed polygonal line~$l_2$ by an elementary cancellation, then $w(l_1,O) = w(l_2,O)$}.

    \begin{figure}[h]\centering
    \compile{
    \includegraphics[scale=0.3]{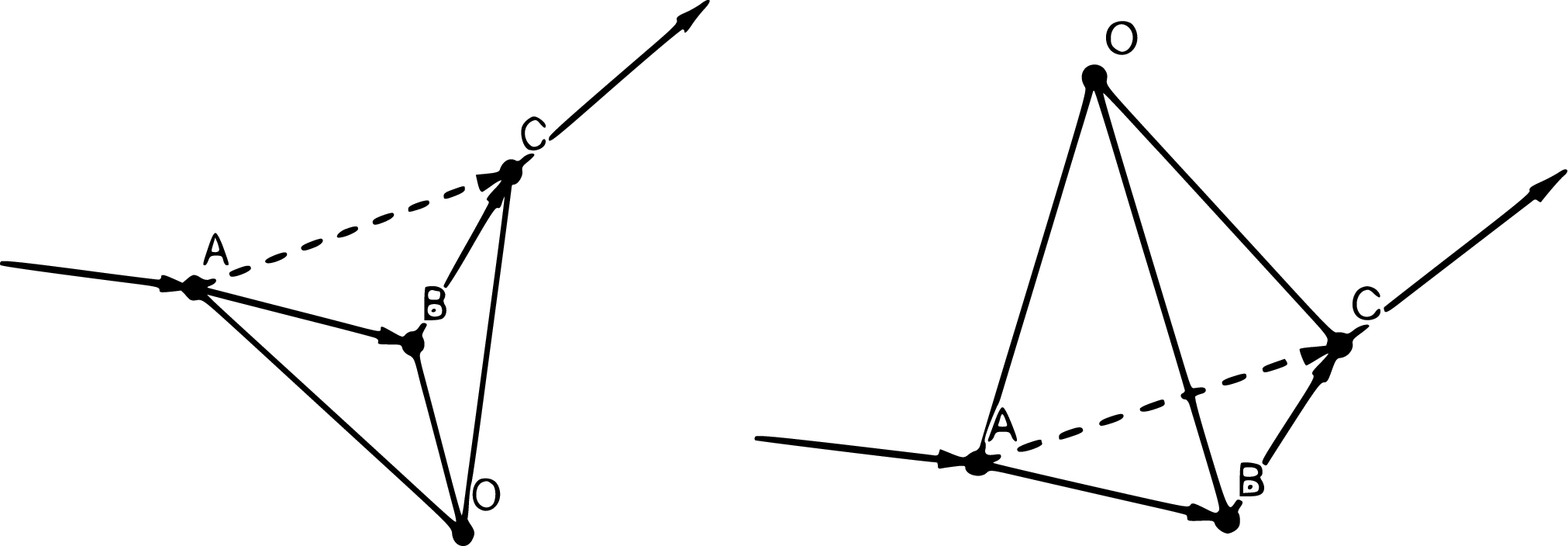}
    }
    \caption{An elementary cancellation of the vertex $B$}
    \label{f:wind}
    \end{figure}

    Let us prove this.
    Let $A,B,C$ be consecutive points of $l_2$ such that `cancellation' of $B$ gives a closed polygonal line~$l_1$.
    Then 
    $$2\pi(w(l_1,O)-w(l_2,O)) = \angle AOB+\angle BOC-\angle AOC = 0,$$
    where the second equality holds by the second property of oriented angles in \S\ref{s:wind} (see Figure~\ref{f:wind}).

    (b) The `only if' part follows from (a).
    To prove the `if' part, it is sufficient to prove that
    \textit{if a closed polygonal line has more than one vertex, and the winding number is zero, then an elementary cancellation can be applied}.

    Let us prove this.
    If some two consecutive vertices in the closed polygonal line coincide, one of them can be removed by an elementary cancellation.
    Now assume that every two consecutive vertices are distinct.
    Since the sum of the angles is zero, there are two adjacent angles $\angle AOB$ and~$\angle BOC$ of different signs.
    Then points~$A$ and~$C$ lie in the same half-plane with respect to the line~$BO$.
    Hence, the triangle~$ABC$ does not contain the point~$O$.
    Therefore, we can make an elementary cancellation of the vertex~$B$.

\begin{figure}[h]\centering
\compile{
\includegraphics[scale=.3]{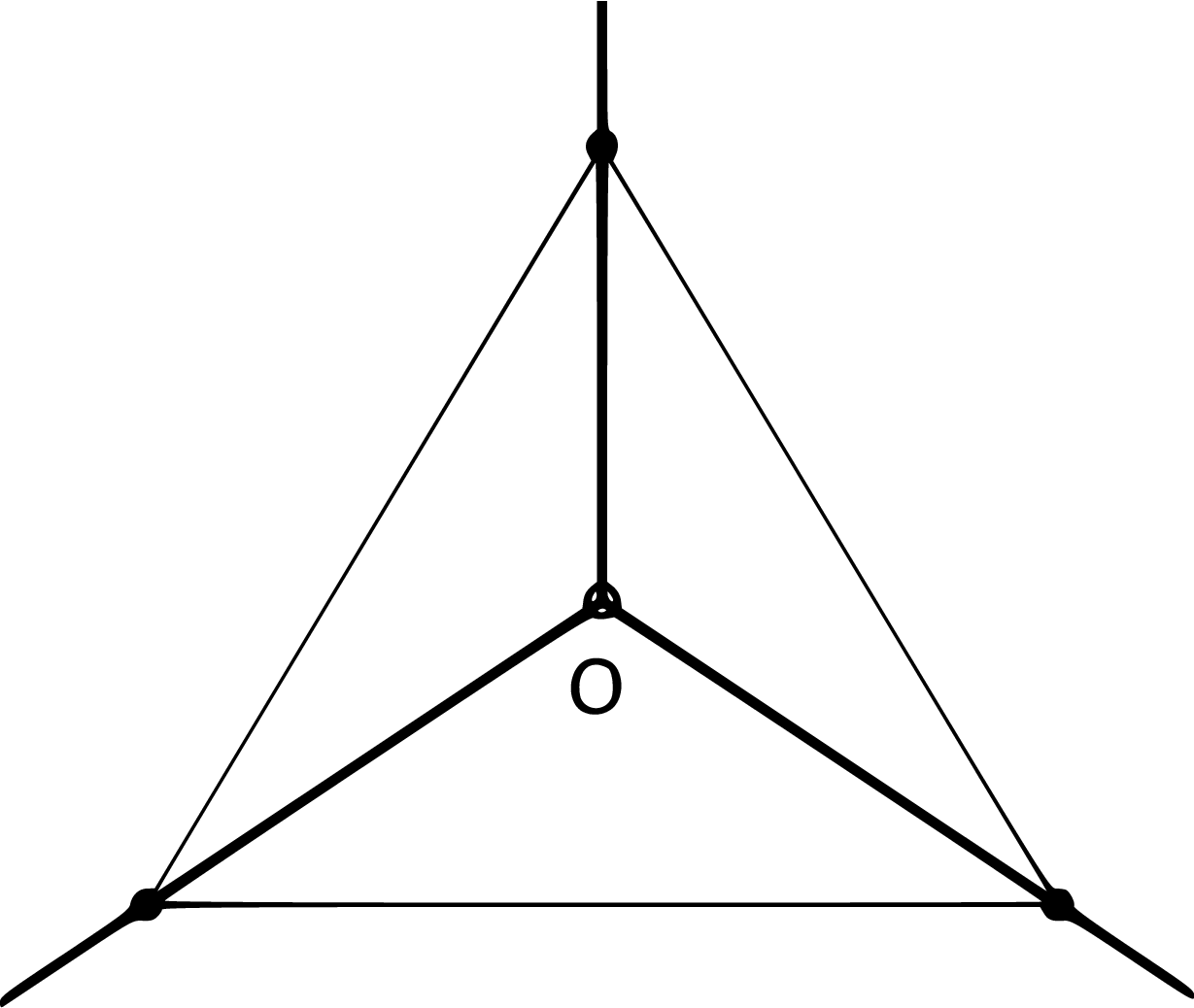}
}
\caption{Additional rays and `dual' triangle in the plane minus a point}
\label{f:compl1}
\end{figure}

    (c) It suffices to prove that {\it any closed polygonal line is homotopic to a one that
    contains only vertices of the triangle $\Delta$} (Figure~\ref{f:compl1}).

    Let us prove this. 
    Denote the vertices of $\Delta$ by $B_1, B_2, B_3$.

    First, in the paragraph after (*) we prove that any closed polygonal line $l$ is homotopic to a closed polygonal line $l'$ such that
  
    (*)  \emph{$l'$ contains at least one vertex of $\Delta$, and for any three consecutive vertices $X, B, Y$ of $l'$ either $B$ is a vertex of $\Delta$, or 
    points $X, B, Y$
    lie in $\angle B_iOB_j$ for some $i,j$.} 

    A pair $\{ XY,\ B \}$ of a segment $XY$ of $l$, and a vertex $B$ of $\Delta$ is called a \emph{bad pair for $l$} if $B \notin \{X, Y\}$ and $B$ lies in $\angle XOY$.
    We may assume that the property (*) does not hold for $l$. 
    We may assume that $l$ contains a vertex of $\Delta$ (by making several 
    inverse operations to an elementary cancellation).   
    Then there is a bad pair $\{ XY,\ B \}$ for $l$.
    Add $B$ to $l$ between $X$ and $Y$ by inverse operation to an elementary cancellation. The number of bad pairs for the obtained closed polygonal line is smaller than for $l$.
    By a finite number of such operations, we can obtain a closed polygonal line $l'$ homotopic to $l$ and satisfying (*).  

    Let $l'$ be a closed polygonal line satisfying (*).
    Then any vertex of $l'$ distinct from $B_1,B_2,B_3$ can be removed by an elementary cancellation. 
    The obtained closed polygonal line also has the property (*). 
    Hence by a finite number of such elementary cancellations, we obtain a closed polygonal line containing only vertices of the triangle $\Delta$
\end{proof}

\section{‘Resolution’ of the Poincar\'e paradox}\label{s:poipar}

The winding numbers (around $P$ and around $Q$) of the example to Problem~\ref{veplho-borr}.a given in Figure~\ref{f:poincare} are obviously zeros. 
Being not homotopic in Problem~\ref{veplho-borr}.a follows from 
Problem~\ref{p:invar-link} (or 
Problem~\ref{p:invar-2eco}).
The idea is to use one of the following homotopy invariants.

\begin{figure}[h]
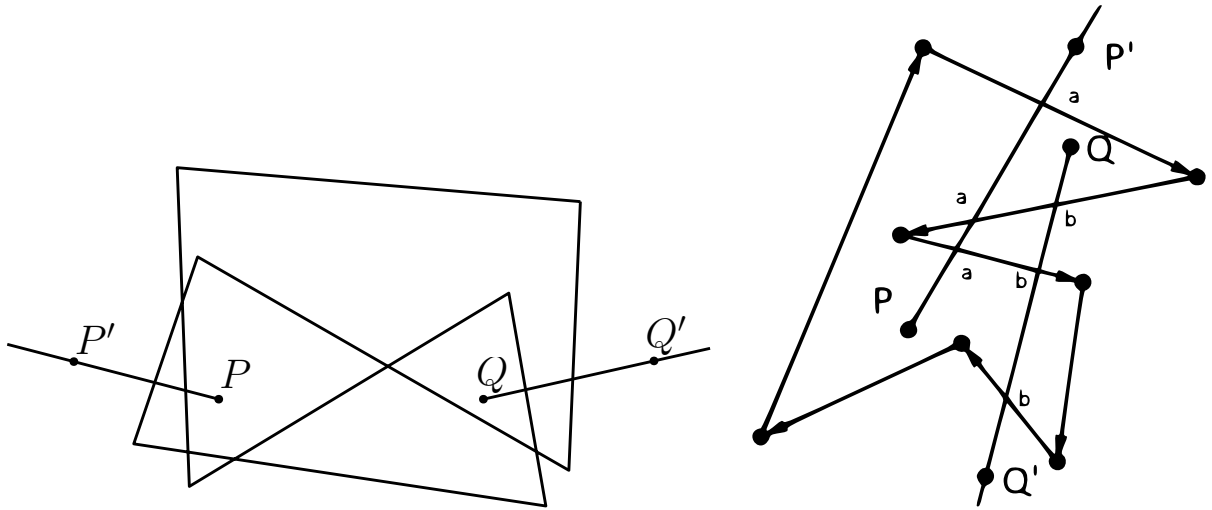
\centering
\compile{
\includegraphics[scale=0.9]{abab.eps}
\quad
\includegraphics[scale=0.4]{cyclic-word.eps}
}
\caption{Cyclic Poincar\'e words modulo 2: $abab$ (left), $abaabb$ (right)}
\label{f:cycl}
\end{figure}

In \S\S\ref{s:poipar}--\ref{s:minus2}

$\bullet$ $l$ is a closed polygonal line in the plane minus the points $P,Q$; 

$\bullet$ $PP'$ and $QQ'$ are disjoint rays $PP'$ and $QQ'$ in general position w.r.t. $l$.

Moving along $l$ write out the letter $a$ (letter $b$) at every intersection with the ray $PP'$ (the ray $QQ'$). 
The resulting oriented cyclic word is called \emph{cyclic Poincar\'e word modulo 2} (see Figure~\ref{f:cycl}).

\begin{pr}\label{p:poin2}
(a) Any cyclic word in letters $a$ and $b$ is the cyclic Poincar\'e word modulo 2 of some closed polygonal line (for some rays $PP'$, $QQ'$).

(b)* The parity of the winding number of a closed polygonal line $l$ around $P$ (around $Q$) equals the parity of the number of letters $a$ (letters $b$) in the cyclic Poincar\'e word modulo 2 of $l$ (for any rays $PP'$, $QQ'$).
\end{pr}

\begin{figure}[h]\centering
\compile{
\includegraphics[scale=1.]{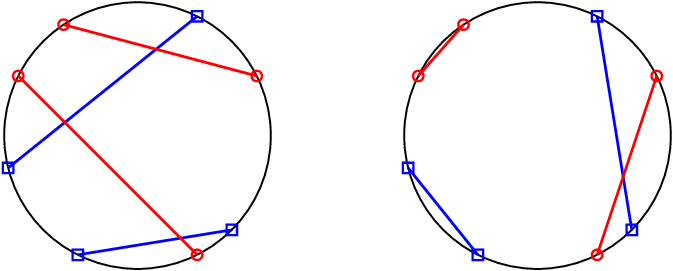}
}
\caption{Two splittings into pairs of four red points (small circles) and four blue points (small squares) on a circle}
\label{f:circles}
\end{figure}

Assume that there are an even number of red points, and an even number of blue points on a circle, and all these points are pairwise distinct. 
Split the red points into pairs, and the blue points into pairs.
Join the points in every red (blue) pair by a red (blue) chord (see Figure~\ref{f:circles}).
The sets of red and blue points are said to be \emph{linked}, if the number of pairs of intersecting red and blue chords is odd (cf.~\cite[Assertion 4.9.5.a]{Sk}).

\begin{pr}\label{p:mil-dis} The property of being linked does not depend on the choice of partitions into pairs (see Figure~\ref{f:circles}).
\end{pr}

A closed polygonal line $l$ in the plane minus points $P,Q$ is said to be \emph{interesting} if $w(l,P)$ and $w(l,Q)$ are even.
By 
Problem~\ref{veplho-poly}.a, the property of being interesting does not change under replacing a closed polygonal line $l$ by a homotopic one.
By 
Problem~\ref{p:poin2}.b, for an interesting closed polygonal line $l$ the number of letters $a$ (letters $b$) in the cyclic Poincar\'e word modulo 2 of $l$ is even.
An interesting closed polygonal line $l$ is said to be \emph{linked} if the letters $a$ and $b$ in the cyclic Poincar\'e word modulo 2 of $l$ are linked.

\begin{pr}\label{p:invar-link} 
For an interesting closed polygonal line, the property of being linked does not change under

(a) changing the rays $PP'$ and $QQ'$;

(b) replacing a closed polygonal line $l$ by a homotopic one.
\end{pr}

\emph{Hint.} Item (b) follows from item (a). 

In 
Problem~\ref{p:invar-link}.b, the cyclic Poincar\'e word modulo 2 may not be defined for some intermediate polygonal lines of homotopy. 
Therefore, for a direct proof of 
Problem~\ref{p:invar-link}.b the analogue of 
Problem~\ref{p:apprgen} for the plane minus two points is useful. 

\section{A homotopy invariant in the plane minus two points}\label{s:poiword}

\begin{pr}\label{p:link-noncomp} There is $l$ such that $w(l, P) = w(l, Q) = 0$,\ $l$ is not linked, but $l$ is not homotopic to a one-point closed polygonal line.
\end{pr}

Consider the set of all (finite) cyclic words (including the empty word) in letters $a,b$.
For such a word, an \emph{elementary cancellation} is the replacement of any of the subwords $aa,bb$ by the empty subword. 
Such a word is said to be \emph{economical} if it has no subwords $aa$ and $bb$. 

\begin{pr}\label{p:wd-2eco} A cyclic word in letters $a,b$ yields by elementary cancellations the unique economical word. 
\end{pr}

An \emph{economical form} (or normal form) of a cyclic word $w$ in letters $a,b$ is the economical word obtained from $w$ by elementary cancellations. 
Denote by $E_{2,2,c}$ the set of all economical words.
Denote by $e_2(l)\in E_{2,2,c}$ the economical form of the cyclic Poincar\'e word modulo 2 of $l$. 
A priori $e_2(l)$ depends on rays $PP'$ and $QQ'$.

\begin{pr}[cf. Problem~\ref{p:invar-link}]\label{p:invar-2eco} 
The word $e_2(l)$ does not change under 

(a) changing the rays $PP'$ and $QQ'$;  \quad 

(b) replacing a closed polygonal line $l$ by a homotopic one.  
\end{pr}

\begin{pr}\label{p:ncomp-2eco}  
There is $l$ such that $e_2(l)$ is the empty word, but $l$ is not homotopic to a one-point closed polygonal line.
\end{pr}

Two cyclic words in letters $a,b$ are said to be \textit{equivalent} if they can be joined by a sequence of cyclic words, in which one of any two consecutive cyclic words can be obtained from the other by an elementary cancellation. 
Denote by $F_{2,2,c}$ the set of all equivalence classes.

\begin{pr}\label{p:word2-example}* (a) The word $abab$ is not equivalent to the empty word.

(b) Are there two non-equivalent cyclic words, both consisting of an odd number of $a$'s, and an even number of $b$'s?

(c) No two distinct economical words are equivalent.

(Hint: use 
Problem~\ref{p:wd-2eco}.)

(d) (Riddle) Construct a `natural' 1--1 correspondence $E_{2,2,c}\to F_{2,2,c}$.

(Hint: use item (c).)

(e) (Riddle) Describe $E_{2,2,c}$ or, equivalently, $F_{2,2,c}$.
\end{pr}

\begin{pr}[cf. Problem~\ref{p:invar-2eco}]\label{p:invar-2eq}* 
The equivalence class in $F_{2,2,c}$ of the cyclic Poincar\'e word modulo 2 does not change under 

(a) changing the rays $PP'$ and $QQ'$;  \quad 

(b) replacing a closed polygonal line $l$ by a homotopic one.
\end{pr}

\section{Homotopy classification in the plane minus two points}\label{s:minus2}

The above homotopy invariants are not complete (see Problems~\ref{veplho-borr},~\ref{p:3points}, \ref{p:link-noncomp},~\ref{p:ncomp-2eco}). 
Now we are ready to introduce a complete invariant.

In this section, closed polygonal lines and their homotopy are considered \emph{in the plane minus two points $P$ and $Q$}.

\begin{figure}[h]\centering
\compile{
\includegraphics[scale=.75]{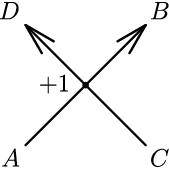}\quad \includegraphics[scale=.8]{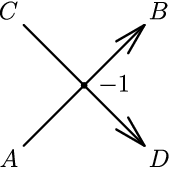}
}
\caption{The sign of intersection point}
\label{f:sign}
\end{figure}

Let $A,B,C,D$ be points in the plane, of which no three belong to a line. 
The {\bf sign} of the intersection point of oriented segments $\overrightarrow{AB}$ and $\overrightarrow{CD}$ is $+1$, if $ABC$ is oriented clockwise, and is $-1$ otherwise (see Figure~\ref{f:sign}).

\begin{pr}[cf. Problem~\ref{p:ray}]\label{p:sign-ray}* Take a closed polygonal line $m$ in the plane minus a point $O$, and a ray $OR$ in general position w.r.t.~$m$. 
Then $w(m,O)$ equals the sum over the oriented segments $x$ of $m$ of the signs of the intersection points of $x$ and $OR$.
\end{pr}

\begin{figure}[H]\centering
\compile{
\includegraphics[scale=0.4]{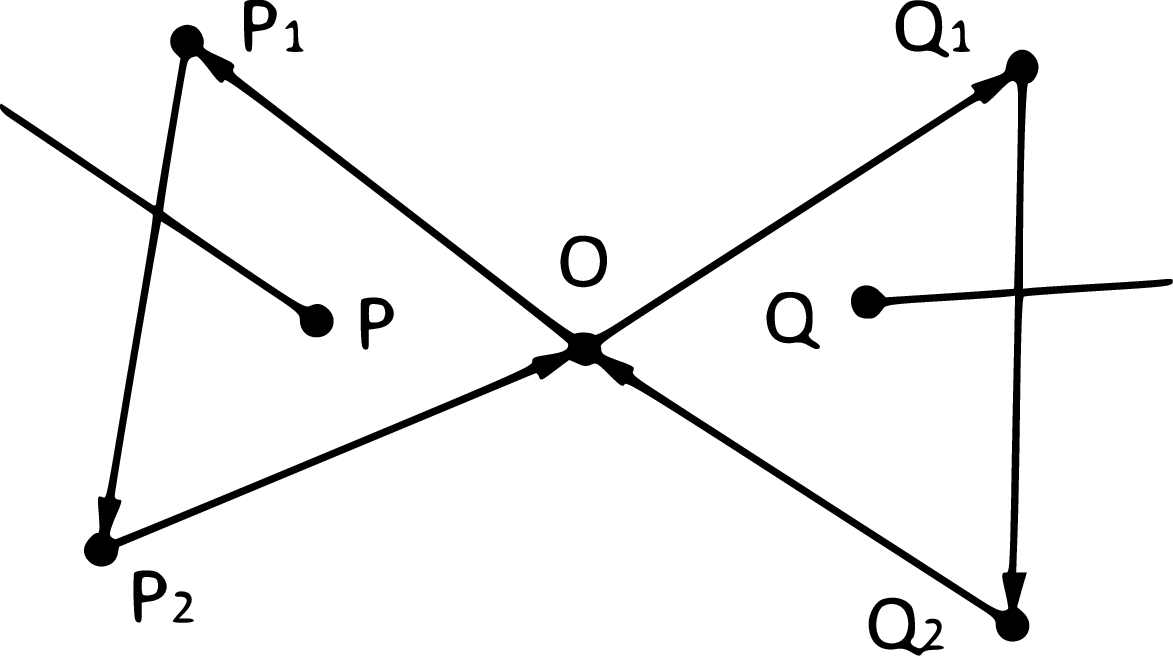}
}
\caption{Figure `eight' and general position rays}
\label{f:eight}
\end{figure}

\begin{pr}\label{p:cycpoi} (a)~(Riddle) Analogously to \S\ref{s:poipar} define for a closed polygonal line the \emph{cyclic Poincar\'e word} in letters $a,b,a^{-1},b^{-1}$.

(Hint: use the definition of the sign.) 

(b) (cf. Problem~\ref{p:poin2}.a) Any cyclic word in letters $a,b,a^{-1},b^{-1}$ is the cyclic Poincar\'e word of some closed polygonal line.

(Hint: use Figure~\ref{f:eight}.)

\end{pr}

Consider the set of all (finite) cyclic words (including the empty word) in letters $a,a^{-1},b,b^{-1}$.
For such a word an \emph{elementary cancellation} is the replacement of any of the subwords $aa^{-1},a^{-1}a,bb^{-1},b^{-1}b$ by the empty subword.
Such a word is said to be \emph{economical} if it has neither of these subwords.

\begin{pr}\label{p:wd-eco} A cyclic word in letters $a,b, a^{-1}, b^{-1}$ yields by elementary cancellations the unique economical word.
\end{pr}

An \textit{economical form} (or normal form) of a cyclic word $w$ in letters $a,b, a^{-1}, b^{-1}$ is the economical word obtained from $w$ by elementary cancellations.
Denote by $E_{2, c}$ the set of all economical words. 
Denote by $e(l)\in E_{2, c}$ the economical form of the cyclic Poincar\'e word of $l$. 
A priori $e(l)$ depends on rays $PP'$ and $QQ'$.

\begin{pr}[cf. Problem~\ref{p:invar-2eco}]\label{p:invar-eco}
The word $e(l)$ does not change under

(a) changing the rays $PP'$ and $QQ'$;  \quad 

(b) replacing a closed polygonal line $l$ by a homotopic one.
\end{pr}

\begin{pr}[cf. Problem~\ref{p:ncomp-2eco}]\label{p:comp-eco}
If $e(l_1)=e(l_2)$ for closed polygonal lines $l_1$ and $l_2$, then $l_1$ and $l_2$ are homotopic. 
\end{pr}

\begin{figure}[h]\centering
\compile{
\includegraphics[scale=.3]{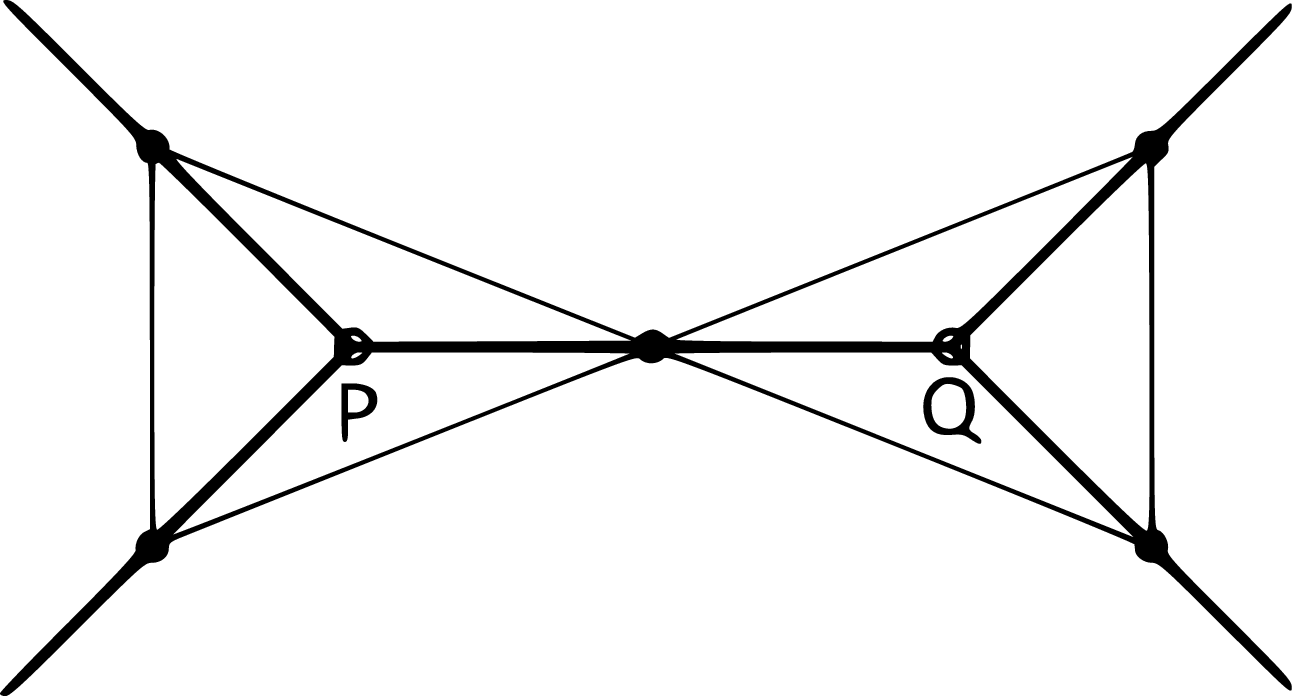}
}
\caption{Additional rays and dual figure `eight' in the plane minus two points}
\label{f:compl2}
\end{figure}

\emph{Hint.} See Figure~\ref{f:compl2} and Problem~\ref{p:defret}. 

\begin{pr}\label{p:defret} (Riddle) Define the property of being homotopic for oriented cycles in the figure-eight graph. 
Construct a 1--1 correspondence between homotopy classes of closed polygonal lines in the plane minus two points, and homotopy classes of oriented cycles in the figure-eight graph.
\end{pr}

\begin{proposition}\label{pr:bij-2}
    The following map between the set of homotopy classes of closed polygonal lines in the plane minus two points, and $E_{2, c}$ is a 1--1 correspondence.
    The image of a closed polygonal line $l$ is the economical form of the cyclic Poincar\'e word of $l$.

    (Proposition~\ref{pr:bij-2} is proved in Problems~\ref{p:cycpoi}.b,~\ref{p:invar-eco}~and~\ref{p:comp-eco}.)
\end{proposition}

This proposition implies Theorem~\ref{t:hocl-algor} for $n = 2$. The analogue of Proposition~\ref{pr:bij-2} for the figure-eight graph instead of the plane minus two points is correct.

\smallskip
\textbf{Remark.} Analogously to \S\ref{s:poiword}, define an equivalence relation on cyclic words in letters $a,b,a^{-1},b^{-1}$. Denote by $F_{2,c}$ the set of all equivalence classes.
Analogously to 
Problem~\ref{p:word2-example}, there is a `natural' 1--1 correspondence $E_{2, c}\to F_{2, c}$.
\smallskip

\begin{figure}[H]\centering
\compile{
\includegraphics[scale=1.]{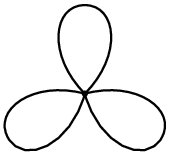} 
}
\caption{The `wedge of 3 cycles' graph}\label{f:wedge}
\end{figure}

\begin{pr}\label{p:bij-n}* (a) The same as in Proposition~\ref{pr:bij-2} replacing the plane minus two points by a plane minus $n$ points.

(b) The same as in Problem~\ref{p:defret} replacing the plane minus two points by a plane minus $n$ points, and the figure-eight graph by the `wedge of $n$ cycles' graph (see Figure~\ref{f:wedge}).

(c) (Riddle) Is the analogue of Problem~\ref{p:invar-eco}.a for the plane minus three points true?
\end{pr}

\emph{Hint.} For a proof of (a), the analogue of 
Problem~\ref{p:apprgen} for the plane minus $n$ points is useful.

\section{Multiplication of closed polygonal lines}\label{s:mulpol}

A \textbf{based closed} oriented \textbf{polygonal line} $A_1\ldots A_m$ in the plane is a sequence (i.e. ordered set)  $(A_1,\ldots,A_m)$ of points in the plane (the points need not be distinct). 
Below, the words `oriented' and `in the plane' are omitted. 
The point $A_1$ is called a \emph{basepoint}.

The \textbf{inverse} to a based closed polygonal line $l = A_1 A_2 \ldots A_m$ is the based closed polygonal line $l^{-1} := A_1 A_m \ldots A_2$.

In this section, $l_1,l_2$ are based closed polygonal lines with a common basepoint $X$.

\begin{figure}[H]\centering
\compile{
\includegraphics[scale=0.3]{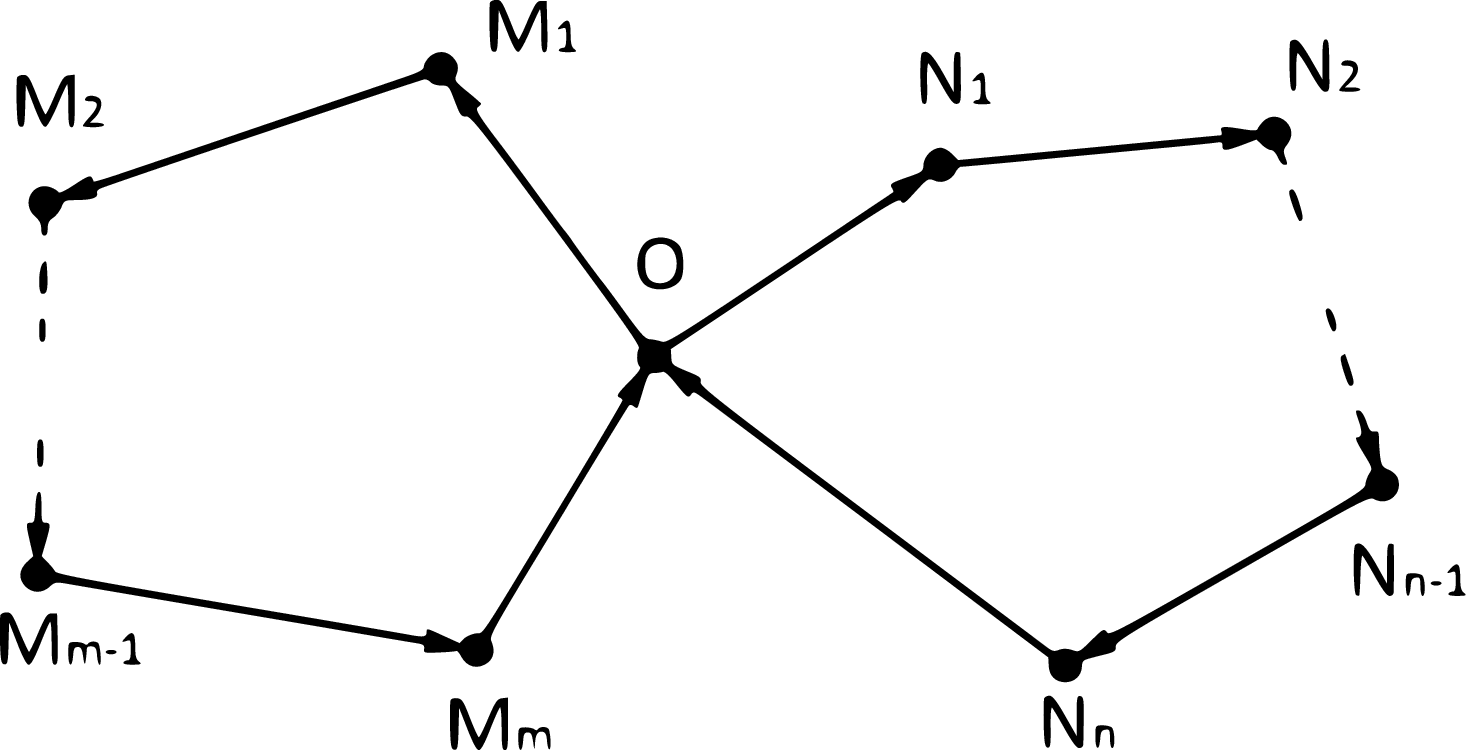}
}
\caption{The product of based closed polygonal lines $l_1 = OM_1\ldots M_m$ and $l_2 = ON_1\ldots N_n$}
\label{f:konk}
\end{figure}

The \textbf{product} (concatenation, joining) of based closed polygonal lines $l_1 = XM_1\ldots M_m$ and $l_2 = XN_1\ldots N_n$ having a common basepoint (see Figure~\ref{f:konk}) is the based closed polygonal line
$$l_1l_2 := XM_1 \ldots M_m X N_1 \ldots N_n.$$ 

\smallskip
\textbf{Remark.} In the plane minus two points $P,Q$ take based closed polygonal lines $a$ and $b$ (triangles in Figure~\ref{f:eight}) with a common base point, and whose convex hulls intersect $\{P,Q\}$ at $P$ and at $Q$, respectively. 
In this remark, we shorten `the closed polygonal line obtained from a based closed polygonal line $x$ by forgetting basepoint' to `a closed polygonal line $x$'.

An example to Problem~\ref{veplho-borr}.a is the closed polygonal line $[a,b] := aba^{-1}b^{-1}$. 
An example to Problem~\ref{veplho-borr}.b is the closed polygonal line $[a,b]c[a,b]^{-1}c^{-1}$ (make this rigorous by yourself).

A version of Poincar\'e paradox is as follows. 
\emph{The closed polygonal lines $aba^{-1}$ and $b$ are homotopic, but the closed polygonal lines $aba^{-1}b^{-1}$ and $bb^{-1}$ are not homotopic ($bb^{-1}$ is homotopic to a one-point closed polygonal line).}

\begin{pr}[additivity]\label{p:3rays} The winding number of the product $l_1l_2$ around a point $O$ outside $l_1,l_2$ is equal to the sum of the winding numbers of $l_1$ and $l_2$ around $O$.
\end{pr}

\begin{pr}\label{p:ngroup} The multiplication of based closed polygonal lines with a common basepoint $X$

(a) is associative, but has no identity element; 

(b) (Riddle) does not generate a well-defined product on homotopy classes.
\end{pr}

Based closed polygonal lines with a common basepoint are said to be \emph{based homotopic} if they are `homotopic preserving basepoint' (give a rigorous definition yourself).

\begin{pr}\label{p:based-nonbased} (a) If $l_1$, $l_2$ are homotopic in the plane minus a point $O$, then they are also based homotopic in the plane minus the point $O$.

(b) There are based closed polygonal lines in the plane minus two points that are homotopic but are not based homotopic.
\end{pr}

\begin{pr}\label{p:group}
(a) The multiplication of based homotopy classes of based closed polygonal lines is well-defined, is associative, and has an identity element.

(b) For this multiplication, every class has an inverse.
\end{pr}

So a natural group structure exists not on the set we are interested in (the set of closed polygonal lines, up to homotopy) but on the set that currently seems less natural (the set of based closed polygonal lines, up to based homotopy). 
Indeed, to define a product, we need \emph{based} closed polygonal lines; 
to obtain a well-defined operation on homotopy classes we need \emph{based} homotopy (Problems \ref{p:ngroup}.b and \ref{p:group}). 

Below, based closed polygonal lines and their based homotopy are considered \emph{in the plane minus two points $P$ and $Q$}.

The \emph{Poincar\'e word} of a based closed polygonal line $x$ is defined analogously to the cyclic Poincar\'e word, only the movement along $x$ begins and ends at the basepoint. 

\begin{pr}[cf. Problem~\ref{p:cycpoi}.b]\label{p:surj-2-based}* Any word in letters $a,b,a^{-1},b^{-1}$ is the Poincar\'e word of some based closed polygonal line.
\end{pr}

Analogously to~\S\ref{s:minus2}, one can define the economical form of a word in letters $a,b,a^{-1},b^{-1}$. 
The economical form of the Poincar\'e word of a based closed polygonal line is a based homotopy invariant (analogously to Problems~\ref{p:invar-2eco},~\ref{p:invar-eco}). 
Below, we present an alternative definition of a based homotopy invariant, in terms of an equivalence relation on the set of all words in letters $a,b,a^{-1},b^{-1}$ (analogous to Problem~\ref{p:invar-2eq}).

Consider the set of all (finite) words (including the empty word) in letters $a,b,a^{-1}, b^{-1}$.
For such a word, an \emph{elementary cancellation} is the replacement of any of the subwords $aa^{-1}, a^{-1}a, bb^{-1}, b^{-1}b$ by the empty subword. 
Two such words are said to be \textit{equivalent}, if they can be joined by a sequence of words, in which one of any two consecutive words can be obtained from the other by an elementary cancellation. 
Denote by $F_2$ the set of all equivalence classes (standard notation: $\left<a,b\right>$).

\begin{pr}[cf. Problem~\ref{p:invar-2eq}]\label{p:invar-eq-based} *
(a) Is it true that the equivalence class of the Poincar\'e word of a based closed polygonal line does not change under changing the rays $PP'$ and $QQ'$?

(b) The equivalence class of the Poincar\'e word does not change under replacing a based closed polygonal line $l$ by a based homotopic one. 
\end{pr}

\begin{pr}\label{p:group2}*
    The multiplication of words generates an operation on $F_{2}$, and turns $F_{2}$ into a group (this group is called the \emph{free group on two generators}).
\end{pr}

\begin{pr}[cf. Problem~\ref{p:comp-eco}]\label{p:comp-eq-based} *
    If the Poincaré words of two based closed polygonal lines $l_1, l_2$ are equivalent, then $l_1$ and $l_2$ are based homotopic.
\end{pr}

\textbf{Remark.} Notice that $F_{2,c}$ (see definition in the remark after Proposition~\ref{pr:bij-2}) is the set of \emph{conjugacy} classes of the group $F_2$.

\begin{proposition}[cf. Proposition~\ref{pr:bij-2}]\label{pr:bij-2-based} 
The following map between the set of based homotopy classes of based closed polygonal lines in the plane minus two points, and $F_2$ is an isomorphism of groups.
The image of a based closed polygonal line $x$ is the equivalence class of the Poincar\'e word of $x$.

(This is proved in Problems~\ref{p:surj-2-based}--\ref{p:comp-eq-based}.)
\end{proposition}

\begin{theorem}[riddle]\label{hompi-alg} (a)~For any graph $G$ with a basepoint, there is an algorithm recognizing whether based oriented cycles in~$G$ are based homotopic.

(b) For any graph $G$, there is an algorithm recognizing whether oriented cycles in~$G$ are homotopic.

(c) For any two graphs, there is an algorithm recognizing whether \emph{simplicial mappings} between graphs are homotopic. 
(See definitions in \cite[\S\,9 `Homotopy classification of maps']{Sk}, \cite{Sk20e}.)
\end{theorem}

Analogues of this theorem for $2$-hypergraphs are false! 
This follows from the algorithmic unsolvability of the triviality problem in some group defined by a finite number of generators and relations. 
See~\cite[Theorem 14.3.1]{Sk20} and a generalization in~\cite{Sk20e}.


%% file: algor_style.tex
Основные идеи представлены на <<олимпиадных>> примерах: на простейших
частных случаях, свободных от технических деталей, и со сведением к необходимому минимуму научного языка.
Благодаря этому и текст становится доступным для начинающих, и удается быстро добраться до интересных сложных и важных результатов, методов и теоретических идей.

Таким образом, многие из приведенных идей работают для более общих случаев.
Мы не тратим время читателя на несложные обобщения.
Их легко придумать или найти в литературе.
Трудно
именно \emph{применить} общую теорию для ярких результатов, сформулированных вне этой теории, когда неизвестно заранее, какую теорию нужно применять (и вообще, можно ли хоть какую-то теорию применить).

\algor{

\emph{`Listeners are prepared to accept unstated (but hinted) generalizations much more than they are able, on the spur of the moment, to decode a precisely stated abstraction and to re-invent the special cases that motivated it in the first place.'}\footnote{Слушатели в гораздо большей степени готовы воспринимать не сформулированные (но упомянутые) обобщения, чем способны <<навскидку>> раскодировать явно сформулированную абстракцию и воспроизвести те конкретные примеры, которые изначально ее мотивировали.} \cite{Ha74}

\emph{`The preceding theorem generalizes to many situations. In fact, there may not actually be an optimal generalization in the sense that no matter what generalization is given, someone could produce a more general one. One of the advantages of being an understander rather than a theorem quoter is that one may be able to obtain approaches to a wide variety of theorems some of which may not even have been formulated yet.'}\footnote{Предшествующая теорема имеет многочисленные обобщения. В действительности, «наилучшего» обобщения может и не быть, в том смысле, что для каждого обобщения кто-то сможет найти еще более общее. Одним из преимуществ читателя, ориентированного на понимание сути, а не просто на цитирование теорем, является способность получать подходы к широкому спектру теорем, некоторые из которых могут быть даже еще не сформулированы.} \cite{Bi83}

}

Более простой материал приводится, чтобы сделать естественным и доступным более сложный.
Попытка начинать с более простого (например, с частных случаев) повышает самостоятельность --- а, значит, глубину и надежность --- освоения материала.
Проще самому доказать частный случай, самому продумать переход от частного к общему, чем  самому сразу доказать общий случай.
Самостоятельно придуманное надежнее запоминается и легче модифицируется.
Кроме того, обычно на частном случае проще отловить и исправить ошибки.
Подробнее см. \cite[\S11.2]{ZSS}.


Как правило, мы приводим \emph{формулировку} красивого или важного утверждения \emph{перед} последовательностью определений и результатов, составляющих его \emph{доказательство}.\algor{\footnote{Часто происходит обратное: формулировки красивых результатов и важных проблем, ради которых была придумана теория, приводятся только {\it после} продолжительного изучения этой теории (или не приводятся совсем). Это развивает представление о математике как науке, изучающей немотивированные понятия и теории. Такое представление принижает ценность математики.}}
В~таких случаях для доказательства утверждения требуется часть дальнейшего материала.
Об этом указано после формулировки утверждения.
Некоторые результаты приведены без доказательств, тогда даются ссылки вместо указаний. 
Если к утверждению (или задаче) не приведено ни доказательство, ни ссылка на него, то оно несложно.
Основные результаты называются <<теоремами>>, менее важные результаты "--- <<утверждениями>>, важные вспомогательные результаты "--- <<леммами>>.
\lktg{Если не оговорено противное, то все это (а также примеры) --- задачи для решения.}

В тексте есть небольшое количество задач (полезно рассматривать и утверждения как задачи).
\algor{Изучение путем решения задач характерно для серьезного обучения математике, см. \cite[\S1.1]{HC19}, \cite[\S1.2]{ZSS} и данные там ссылки.
Оно продолжает древнюю культурную традицию.\footnote{Например, послушники дзенских монастырей обучаются, размышляя над загадками, данными им наставниками.
Впрочем, эти загадки являются скорее наводящими на размышления парадоксами, а не задачами.
См. подробнее \cite{Su}.}
Задачи пронумерованы жирными цифрами, слово <<задача>> опускается.
In a lecture we often start with a question and only then give an answer, e.g.
`Может ли количество ... быть равным 7?'  instead of Assertion \ref{pla-kon}.a.}
Если условие задачи является формулировкой утверждения, то в~задаче требуется это утверждение доказать (и тогда в ссылках мы называем это утверждение утверждением, а не задачей).
\algor{
Краткие указания
к задачам
приводятся сразу, более длинные в конце параграфа.}
\emph{Загадкой} называется не сформулированный четко вопрос; здесь нужно придумать и~четкую формулировку, и~доказательство.
\algor{Если номер задачи помечен звездочкой,
то эта задача посложнее соседних, или менее важна для дальнейшего.}

Определения важных понятий даны \textbf{жирным шрифтом}, чтобы их было проще найти.
\algor{{\it Разделы и задачи, отмеченные звездочкой, а также замечания,} не используются в дальнейшем.}

\lktg{Мы приглашаем школьников представлять свои письменные решения на конференциях школьников, см., например, \cite{ZS25}.} 